\documentclass[10pt]{amsart}
\usepackage{amssymb}

\numberwithin{equation}{section}

\def\cro{{\rm cr }}
\def\fdt{{\rm fdt}}

\def\b{{\beta}}

\def\N{{\Bbb N}}

\def\R{{\Bbb R}}

\def\E{{\Bbb E}}

\def\bD{{\bf \Gamma}}
\newcommand{\reals}{{\Bbb{R}}}

\def\IJ{{\bf I}}
\def\wIJ{{\widehat\IJ}}

\def\ZZ {\hbox{\it Z\hskip -4pt Z}}

\def\CC{{\mathcal C}}

\def\CC{{\mathcal R}}
\def\CC{{\mathcal C}}

\def\Aa{{\mathcal A}}
\def\Ba{{\mathcal B}}

\def\Sa{{\mathcal S}}

\def\R{{\Bbb R}}

\def\N{{\Bbb N}}
\def\E{{\Bbb E}}

\def\prf{\noindent {\bf Proof:\ }}

\def\nn{\noindent}

\def\hh{\widehat}

\def\b{\beta}
\def\d{\delta}

\def\s{\sigma}

\def\NC{{\mbox{NC}}}
\def\D{\Delta}

\def\L{\Lambda}

\def\part{\partial}

\def\bu{\bullet}

\def\ts{\times}

\def\ra{\rightarrow}

\def\tilde{\widetilde}
\def\hh{\widehat}
\def\npsi{\hh\psi}
\def\phsi{\varphi}

\cleardoublepage
\newtheorem{prop}{Proposition}[section]

\newtheorem{lem}[prop]{Lemma}

\newtheorem{theo}[prop]{Theorem}

\begin{document}              

\title[Limiting dynamics at high temperature]
{Limiting dynamics for spherical models of spin glasses \\
at high temperature}
\author{Amir Dembo}
\address{Department of Statistics and Department of Mathematics\\
Stanford University\\ Stanford, CA 94305.}
\email{amir@math.stanford.edu}

\author{Alice Guionnet}

\address{UMPA, Ecole Normale Superieure de Lyon\\
46 all\'ee d'Italie\\
69364 Lyon Cedex 07, France}
\email{Alice.Guionnet@umpa.ens-lyon.fr}

\author{ Christian Mazza}
\address{D\'epartement de math\'ematiques, Universit\'e P\'erolles, Chemin du Mus\'ee 23, CH-1700 Fribourg, Suisse}
\thanks{\noindent Research partially supported by NSF
grants \#DMS-0406042, \#DMS-FRG-0244323. Supported in parts by the Swiss National Foundation.
\newline
{\bf AMS (2000) Subject Classification:}
{Primary: 82C44 Secondary:  82C31, 60H10, 60F15, 60K35}
}

\begin{abstract}
We analyze the coupled non-linear integro-differential 
equations whose solution is the thermodynamical limit of 
the empirical correlation and response functions 
in the Langevin dynamics for spherical $p$-spin 
disordered mean-field models. We provide a mathematically
rigorous derivation of their FDT solution (for   
the high temperature regime) and of certain key 
properties of this solution, which are in agreement with 
earlier derivations based on physical grounds. 
\end{abstract}

\maketitle

\section{Introduction\label{s.introduction}}

The complex long time behavior 
predicted for the thermodynamical limits of a wide 
class of Markovian dynamics with random interactions,
is among the fascinating aspects of out of equilibrium 
statistical physics (for a good survey on phenomena 
such as aging, memory, rejuvenation, and violation of the 
Fluctuation-Dissipation Theorem (FDT), see \cite{BKM,LesHouches}). 
This work is concerned with the long time behavior
of a complex system composed of $N$ Langevin particles
$x_t=(x_t^i)_{1\le i\le N}\in\R^N$, each evolving in
$\R$ and interacting with the others through a random
potential. More precisely, one considers a diffusion of the form
\begin{equation}\label{diffusion}
 dx_t=-f'(\vert\vert x_t\vert\vert^2/N)x_t dt
  -\beta\nabla H_J(x_t) dt +dB_t,
  \end{equation}
  where $B_t$ is a $N$-dimensional Brownian motion, $\vert \vert x \vert
\vert$ denotes the Euclidean norm of $x \in \R^N$ and
$f$ is a convex function. Such models are called
{\it spherical} (for well chosen $f$, $x_t$ is restricted
to stay on a sphere). The mixed $p$-spin, $p \leq m$,
potential $H_J:\ \R^N\longrightarrow \R$ is given by
\begin{equation}\label{potential}
  H_J(x)=\sum_{p=2}^m \frac{a_p}{p!}\sum_ {1\le i_1\le \cdots\le
i_p\le N}J_{i_1\cdots i_p}x^{i_1}\cdots x^{i_p},\quad a_m\neq 0
\end{equation}
where the coupling constants $J_{i_1\cdots i_p}$ are assumed to be
independent centered Gaussian variables. The variance of $J_{i_1\ldots i_p}$
is $c(\{i_1,\ldots,i_p\}) N^{-p+1}$, where
\begin{equation}\label{eq:vardef}
c(\{i_1,\ldots,i_p\})=\prod_k l_k! \,,
\end{equation}
and $(l_1,l_2,\ldots)$ are
the multiplicities of the different elements of the set      
$\{i_1,\ldots,i_p\}$ (for example, $c=1$ when $i_j \neq i_{j'}$ for
any $j \neq j'$, while $c=p!$ when all $i_j$ values are the same).

When $m=2$, one gets the so-called
{\it Sherrington-Kirckpatrick} spherical spin glass,
  which has been studied in details
in \cite{2001}. Given a realization of the coupling
  constants, the dynamics of (\ref{diffusion}) is
invariant for the (random) Gibbs measure
  \begin{equation}\label{Gibbs}
  \mu_N^J(d x)=Z_{J,N}^{-1}
\exp(-N f(\vert\vert x\vert\vert^2/N)-2\beta H_J(x))\prod_{i=1}^N d x^i.
  \end{equation}
Similar random measures have been
extensively studied in mathematics and physics during
the last two decades (see e.g.
   \cite{Talagrand}, for the rigorous analysis of the asymptotics of
the free energy of the measure with a hard spherical constraint, 
corresponding to spins on the sphere $\vert\vert x\vert\vert^2=N$).
  Here, we shall be concerned with the
statistical properties of the dynamics at high temperature. The
  natural quantity of interest is the {\it empirical covariance function}
  \begin{equation}\label{empiricalcovariance}
  C_N(s,t)=\frac{1}{N}\sum_{i=1}^N x_s^i x_t^i,\ \ s\ge t,
  \end{equation}
  in the large $N$ limit, and for large $t$ and $s$. It turns out that the
asymptotic behavior of
  (\ref{empiricalcovariance}) strongly depends on the way $t$ and $s$ tend
to infinity, at least
  at low temperature. This is a trace of {\it aging}: the older it gets,
the longer the system
  will take to forget its age. This innocent looking notion of aging is
related to deep
  mathematical problems (see e.g. the survey in \cite{Alice}), and leads
to interesting mathematical scenarios, like the ultrametric
  property of the covariance function at low temperature (see e.g. \cite{CK}).
  In \cite{2001}, the authors present a detailed analysis of the aging
properties of
  (\ref{empiricalcovariance}) in the special case $m=2$, using integro
differential equations
  involving the almost sure limit $C(s,t)=\lim_{N\to\infty}C_N(s,t)$. 
When $m\ne 2$, closed equations
for $C$ are obtained in \cite{CK} (for
the hard spherical constraint) and rigorously derived in
\cite{BDG2}, where they are called {\it Cugliandolo-Kurchan equations}.
These equations involve also the limit of the
{\it integrated response function}
\begin{equation}\label{integrated}
 {\chi}_N(s,t)=\frac{1}{N}\sum_{i=1}^N x_s^i B_t^i\,.
\end{equation}
According to \cite{BDG2}, fixing $T<\infty$, the random functions
$C_N$ and ${ \chi}_N$ converge uniformly on $[0,T]^2$, almost surely 
and in $L_1$ to non-random functions $C(s,t)$ and
 ${\chi}(s,t)=\int_0^t R(s,u) du$ with
  $R(s,t)=0$ when $t>s$, $R(s,s)\equiv 1$, and, for $s>t$,
 the absolutely continuous
  functions $C$, $R$ and $K(s)=C(s,s)$ are the unique solutions in the space
of bounded, continuous
  functions, of the non linear integro-differential equations
\begin{eqnarray}\label{Eq1}
\partial_s R(s,t)&\!\!\!\! =\!\!\!\!
& - f'(K(s)) R(s,t) + \b^2 \int_t^s
R(u,t) R(s,u) \nu''(C(s,u)) du ,\label{eqR}\\\label{Eq2}
\partial_s C(s,t)&\!\!\!\! = \!\!\!\!
& - f'(K(s)) C(s,t) +
\b^2 \int_0^s C(u,t) R(s,u) \nu''(C(s,u)) du
+ \b^2 \int_0^t \nu'(C(s,u)) R(t,u) du,\label{eqC}\\\label{Eq3}
\partial_s K(s) &\!\!\!\! =\!\!\!\! & -2 f'(K(s)) K(s) + 1 + 2\b^2
\int_0^s \psi(C(s,u)) R(s,u) du , \label{eqZ}
\end{eqnarray}
where $\psi(r)=\nu'(r)+r\nu''(r)$,
\begin{equation}\label{eq:nudef}
\nu(r)=\sum_{p=2}^m \frac{a_p^2}{p!}r^p,\ \ m\ge 2,\ \ a_m\ne 0,
\end{equation}
and the initial condition $K(0)=C(0,0)>0$ is given.
It was shown in Theorem 1.2 of \cite{BDG2} that for
\begin{eqnarray}
f(r)&:=& f_L(r) = L(r-1)^2 + \frac{1}{4k} r^{2k} \,, \qquad k > m/4, \,
k \in \ZZ, \, L \geq 0 \,,
\label{eq:fdef}
\end{eqnarray}
these equations admit a unique solution
$C(s,t)=C(t,s)$, $R(s,t)$ and $K(s)=C(s,s)$
in the space of absolutely continuous
functions on $\{(s,t): 0 \leq t \leq s\}$.\footnote{
More general choice of differentiable
$f(\cdot)$ is allowed in \cite{BDG2}, but we focus here
on the collection given by (\ref{eq:fdef}).}

Whereas \cite{BDG2} rigorously 
derives these equations using stochastic calculus
and concentration inequalities, in the physics literature they are attained 
via the so-called {\it Martin-Siggia-Rose formalism} 
(see e.g. \cite{BCKM} or \cite{CD}). An alternative to the latter is
to expand the stochastic process $x_t$ perturbatively using 
diagrams, as explained for example in \cite{BCKM} and \cite{LesHouches}, 
then average over the disorder the product
$x_s x_t$ to get equations relating the 
covariance and response functions. The generic form 
of the family of diagrams  
one uses in this process indicates that (\ref{Eq1})--(\ref{Eq3})
are equivalent to a class of {\it mode coupling equations}
(see e.g. \cite{BCKM}, \cite{LesHouches} or \cite{KB}).
Mode coupling approximations were developed in physics to study nonlinear 
random dynamical systems
occurring in many contexts like plasma physics, kinetic theory of classical
liquids or glasses (see e.g. \cite{GO2} and the references therein).
This method
considers a perturbative expansion as a series 
indexed by diagrams  containing information on the nonlinearity. 
The series is then renormalized, to produce self-consistent moment
equations, called mode coupling equations.
Kraichnan \cite{KA} \cite{KB} developed such approximations 
(called direct interaction in fluid mechanics) when considering
the solutions to the Navier-Stokes equations in Fourier space under random
initial conditions, and gave a procedure to perform the statistical closure of
the moment equations: usually, there is a cascade of moments meaning that the
time derivative of the second moment involve the third moment and so on. 
One then looks for good approximations leading to self-consistent equations 
(c.f. the texts \cite{MC} and \cite{MY}, for a variety of
statistical closure problems and 
diagrammatic methods in fluid mechanics).
More recently, these approximations were applied in 
the study of super-cooled and strongly interacting liquids (glasses),
producing very accurate quantitative predictions (see  \cite{GO1}, \cite{GO2}). 
The fundamental object of interest in structural glasses is 
a correlation function $\eta(t)$, for which the mode coupling equation is
\begin{equation}\label{eq:mode}
\eta''(t)+\Omega^2 \eta(t)+\nu\eta'(t)+\Omega^2
 \int_0^t k(\eta(t-u))\eta'(u){\rm du} =0\,,
\end{equation}
where $k(r)=\sum_{p=1}^m b_p r^p$, for some constants $\Omega$, $\nu$ and 
$\{b_p, p=1,\ldots,m\}$ (see e.g. \cite{GO1}).
The integro-differential equation (\ref{eq:mode}) is 
similar to what one gets when postulating the FDT ansatz, 
whereby the solution $(R,C)$ of (\ref{Eq1})-(\ref{Eq3})
is translation invariant, 
with $K(s)$ constant and $R$ proportional to the derivative of $C$
(for example, see (\ref{eqRnew}) in the sequel). Of course, the first
task when dealing with (\ref{Eq1})--(\ref{Eq3}) is to show 
the validity of this ansatz, at least for $\beta$ small enough, 
as we do in Theorem \ref{FDT}. Further, the genericity of (\ref{eq:mode})
implies that any relevant information about its 
solution is of interest (and in this context see Proposition 
\ref{studyFDTeq}).

The asymptotic behavior of $C(s,t)$ and $R(s,t)$ for large values of
$t$ and $s$ is difficult to pin down; 
in \cite{CK}, the authors propose various scenarios,
but no complete description of these asymptotics could be given
(see also \cite{BKM} or \cite{LesHouches}). A first regime of interest
is the so-called FDT regime, in which the fluctuation dissipation
theorem of statistical physics is expected to hold. In this regime,
the covariance should be stationary,
that is, for fixed $s-t=\tau$ and $t$ large, the covariance 
$C(s,t)$ 
should be approximated well by some function $C_{\fdt}(\tau)$. 
Further, from \cite{GM} we know that the response function $R(s,t)$
is then well approximated by some function $R_{\fdt}(\tau)$.
In this regime we further expect the FDT relation 
$R_{\fdt}(\tau)=-2 C_{\fdt}'(\tau)$ to hold (the proportionality constant
has to do with the scaling we employed in the definition 
of $R(s,t)$ and not with its physical meaning).
We shall see in this work that the FDT regime holds for small enough
$\beta$ in which case both $C_{\fdt}$ and $R_{\fdt}$ 
decay to zero exponentially fast. The aging regime
is expected to be characterized by covariances and response functions of
the generic form 
$C(s,t)=C_{aging}(h(t)/h(s))$ and $R(s,t)=h'(t) h(s)^{-1} R_{aging}(h(t)/h(s))$
with sub-exponential growth of the monotone function $h(\cdot)$
and a polynomial decay to zero and to appear only for $\b>\b_c$, the
dynamical phase transition point of \eqref{diffusion}.
These different scenarios are examined in \cite{GM}, where 
fixing the asymptotic behavior of $C$ according to the above
choices $C_{\fdt}$ or $C_{aging}$, the authors 
study the solution $\tilde R$ of the equation 
$$
\partial_s \tilde R(s,t)=
 - f'(K(s)) \tilde R(s,t) + \int_t^s
\tilde R(u,t) \tilde R(s,u) k(s,u) du,$$
where $k(s,u)$ stands either for $\beta^2 \nu''(C_{\fdt}(s-u))$
or for  $\beta^2 \nu''(C_{aging}(h(u)/h(s))$.
Let
$$H(s,t)=\exp(\int_t^s f'(K(u)) d u)R(s,t).$$
Then it is easy to check that $H$ solves the equation 
\begin{equation}\label{Kraichnan}
\partial_s H(s,t)=\beta^2\int_t^s H(s,u)H(u,t)\nu''(C(s,u)) d u,
\quad H(t,t)=1\,.
\end{equation}
The main technique of \cite{GM} is a formal solution to \eqref{Kraichnan}
involving non-crossing involutions. We shall also use 
this formula in the sequel
to prove various statements concerning the solution $C$ and $R$
of the complete system of equations 
(for example, their exponential decay at high enough temperature).

We first show that in the limit $L \to \infty$,
the solutions of the equations (\ref{eqR})--(\ref{eqZ})
for our ``soft'' spherical constraint
by the potential $f_L(\cdot)$ coincide with the limiting equations
of \cite{CK} for the hard spherical constraint
(where $C(2\cdot,2\cdot)$ and $R(2\cdot,2\cdot)$ 
are considered for $\nu(r)=r^p/8$).

\begin{prop}\label{prop-sphere}
For any $T<\infty$, the solution $(R_L, K_L, C_L)$
of the system (\ref{eqR})--(\ref{eqZ})
with potential $f_L(\cdot)$ as in (\ref{eq:fdef}) and 
initial condition $K_L (0)=1$, converges
as $L \to \infty$, uniformly in $[0,T]^2$, towards 
the triplet $(R,K,C)$ such that $C(t,t)=K(t)=1$ for all $t \geq 0$,
$R(s,t)=0$ for all $s<t$, and for all $s \geq t$,
\begin{eqnarray}
\partial_s R(s,t)&\!\!\!\! =\!\!\!\!
& - \mu(s) R(s,t) + \b^2 \int_t^s
R(u,t) R(s,u) \nu''(C(s,u)) du ,\label{eqRs}\\
\partial_s C(s,t)&\!\!\!\! = \!\!\!\!
& -\mu(s) C(s,t) +
\b^2 \int_0^s C(u,t) R(s,u) \nu''(C(s,u)) du
+ \b^2 \int_0^t \nu'(C(s,u)) R(t,u) du,\label{eqCs}
\end{eqnarray}
where
\begin{equation}\label{eqZs}
\mu(s)= \frac{1}{2}
\left(1 + 2\b^2
\int_0^s \psi(C(s,u)) R(s,u) du \right).
\end{equation}
Moreover, $C(s,t)=C(t,s)$ is a non-negative definite kernel, with values
in $[0,1]$ and $R(s,t) \geq 0$ is such that $R(t,t)=1$ and 
\begin{equation}\label{eq:rbd-new}
|\int_{t_1}^{t_2} R(s,u) du |^2 \le t_2-t_1 \,,   
\qquad 0 \leq t_1 \leq t_2  \leq s < \infty \,.
\end{equation}
\end{prop}

We prove in Section \ref{expdec}
that for $\beta>0$ sufficiently small, the solution $(R,C)$ of
(\ref{eqRs})--(\ref{eqZs}) decays exponentially fast in $s-t$.
\begin{prop}\label{prop-exp-tails}
There exists $\beta_0>0$ such that for all $\beta < \beta_0$ there
exists $M_\beta < \infty$ and $\eta_\beta>0$ for which the solution 
of (\ref{eqRs})--(\ref{eqZs}) satisfies for all $s \geq t \geq 0$,
\begin{eqnarray}
R(s,t) &\leq&  e^{-\eta_\beta (s-t)} \label{eqRexp} \\
C(s,t) &\leq& M_\beta e^{-\eta_\beta (s-t)} \label{eqCexp} \,.
\end{eqnarray}
\end{prop}

Equipped with Proposition \ref{prop-exp-tails} we prove
in Section \ref{fdtequation} that for some $\beta_1>0$
and each $\beta<\beta_1$, 
when $s-t=\tau$ is fixed and $t\to \infty$ the solution $(R(s,t),C(s,t))$ 
converges to a limiting pair $(R_{\fdt}(\tau),C_{\fdt}(\tau))$ that satisfies 
the FDT relation $R_{\fdt}(\tau)=-2 C_{\fdt}'(\tau)$ for all $\tau \geq 0$.
Upon analyzing in Section \ref{FDTeqsec} the corresponding FDT equations 
\eqref{eqRFDT}--\eqref{eqKFDT}, we establish our main result.
\begin{theo}\label{FDT}
If $\beta<\beta_1 \le 1/(2\sqrt{\nu''(0)})$ then for any $\tau\ge 0$,
$$
\lim_{t\ra\infty}C(\tau+t,t)=C_{\fdt} (\tau),
$$
and
$$
\lim_{t\ra\infty}R(\tau+t,t)=R_{\fdt} (\tau) = - 2 C_{\fdt}'(\tau),
$$
where $C_{\fdt}(\cdot)$ is the unique $[0,1]$-valued,
continuously differentiable solution of the equation
\begin{equation}\label{FDTDb}
 D'(s)=-\int_0^s\phi(D(v)) D'(s-v) dv - b ,\qquad D(0)=1 \,,
\end{equation}
for $b=1/2$ and $\phi(x)=b+2\b^2 \nu'(x)$.
Moreover, both $C_{\fdt}(\cdot)$ and $R_{\fdt}(\cdot)$ decay exponentially
to zero at infinity.
\end{theo}

A key ingredient of the proof of Theorem \ref{FDT} is the analysis of
the equation \eqref{FDTDb}, which is of some independent interest. 
Specifically, in Section \ref{FDTeqsec} we prove the following.
\begin{prop}\label{studyFDTeq}
Suppose $b>0$ and  
$\phi:[0,1] \to \R$ is non-decreasing, of Lipschitz continuous derivative 
such that 
\begin{equation}\label{eq:Dcond}
\sup_{0 \leq x \leq 1} \, \{ \phi(x)(1-x) \} \ge b \,.
\end{equation}
Then, \eqref{FDTDb} has a unique solution in the space of
$[0,1]$-valued continuously differentiable functions. This solution
is twice continuously differentiable, strictly decreasing 
and converges for $s \to \infty$ to 
$$
D_\infty := \sup\{ x\in [0,1]: \phi(x)(1-x)\ge b \} \,.
$$
Further, in case $\phi(\cdot)$ is convex and $\phi(1)>2 \sqrt{b \phi'(1)}$,
the derivative of the solution of \eqref{FDTDb} decays exponentially to zero 
with some positive exponent.
\end{prop}
Setting $b=\phi(D_\infty) (1-D_\infty)$ and considering $x \downarrow D_\infty$ in \eqref{eq:Dcond},
it is not hard to verify that necessarily $\phi(D_\infty) \geq \phi'(D_\infty) (1-D_\infty)$. In
the remark following the proof of Proposition \ref{studyFDTeq} we observe that the condition
\begin{equation}\label{eq:exp-conv}
\phi(D_\infty)>\phi'(D_\infty) (1-D_\infty) \,,
\end{equation} 
is necessary for the exponential convergence of 
$D'(s)$ to zero as $s \to \infty$ when $\phi(\cdot)$ is convex. As \eqref{eq:exp-conv}
is easily seen to be equivalent to $\phi(1)>2 \sqrt{b \phi'(1)}$ when $\phi(x)$
is linear, we speculate that it actually characterizes the exponential convergence
of the solution of \eqref{FDTDb}.

Of course, in case of 
$\phi(\cdot)=b+2\beta^2 \nu'(\cdot)$ and $b=1/2$ the condition 
\eqref{eq:Dcond} holds (try $x=0$). Setting $\b_c \in (0,\infty)$ via
\begin{equation}\label{eqbetac}
\frac{1}{4\beta_c^2}=\sup \{ \nu'(x)(1-x)x^{-1} : x \in (0,1] \} \,,
\end{equation}
it is easy to check that in this case $D_\infty=0$ if 
$\b<\b_c$ whereas $D_\infty>0$ for $\b>\b_c$. 
Further, considering $x \to 0$ in \eqref{eqbetac}
we find that $1/(4\b_c^2) \geq \nu''(0)$, so 
the condition \eqref{eq:exp-conv} then holds for any $\b<\b_c$.
This indicates that though the values of $\beta_0$ and $\beta_1$ 
in our proofs of Proposition \ref{prop-exp-tails} 
and Theorem \ref{FDT} are quite small, 
both should match the predicted 
dynamical phase transition point $\b_c$ of our model.

Indeed, subject to a 
heuristic ansatz, similar to what can be found in the physics literature, 
we show in Section \ref{sec-heur} why the equation \eqref{FDTDb} 
with $b=1/2$ and $\phi(x)=\gamma + 2\b^2\nu'(x)$
describes
the FDT regime at all temperatures, albeit in general with 
$\gamma=\gamma(\b)\neq b$ (see Proposition \ref{ansatz}). 
More precisely, the physics prediction
for the choice of $\gamma$ is as follows.
Let $x^*$ denote the largest value of $x < 1$
attaining the supremum in \eqref{eqbetac}
of the polynomial $h(x)=\nu'(x)(1-x)/x$ (with $h(0)=\nu''(0)$).
Setting $g(x)=\nu''(x) (1-x)^2$ we
have that $g(x^*)=h(x^*)=1/(4 \b_c^2)$ as a
consequence of the optimality condition $h'(x^*)=0$ if $x^*>0$
(while trivially $g(0)=h(0)$). From the continuity of $g(\cdot)$
and the fact that $g(1)=0$ it follows that
\begin{equation}\label{qEA}
q(\b) := \sup\{ x \in [0,1] : 4 \b^2 \nu''(x) (1-x)^2 \geq 1 \}\,,
\end{equation}
is in $[x^*,1)$ for all $\b \geq \b_c$, with $q(\b)>x^*$ as soon
as $\b>\b_c$. Setting $\gamma(\b)=2\b^2 [\nu''(x)(1-x)-\nu'(x)]$ 
for $x=q(\b)$ we find that the 
condition \eqref{eq:Dcond} of Proposition \ref{studyFDTeq} applies 
also for $\b \geq \b_c$. Further, the expected limit $D_\infty$ 
of $C_{\fdt}(\tau)$ as $\tau \to \infty$ is the preceding 
$q(\b)$ which is strictly positive for $\b>\b_c$, indicating the
onset of the aging regime at $\b_c$. 
Another indication of the
onset of aging is the fact that when $\b \geq \b_c$
the exponential convergence of $C_{\fdt}(\cdot)$ and
$R_{\fdt}(\cdot)$ is lost (i.e. the preceding
choice of $\gamma(\b)$ leads to equality in \eqref{eq:exp-conv}).
The physics prediction suggests also that 
$$
\IJ_{\gamma(\b)} =  \b^2 \int_0^1 \psi(C_{aging}(\lambda)) 
R_{aging}(\lambda) d\lambda \geq 0 \,,
$$
where $\IJ_{\gamma(\b)}=\gamma(\b)-\frac{1}{2}+2\b^2 q(\b) \nu'(q(\b))$ 
is strictly positive when $q(\b)>0$, and in particular, 
whenever $\b > \b_c$ (c.f. the remark following Proposition \ref{ansatz}). 
The precise 
nature of the dynamical phase transition at $\b_c$ depends on
whether $q(\b_c)$ is strictly positive (as is the case for 
example whenever $a_2=0$), or not. 
The physics ansatz of one aging regime 
with $R_{aging}(\lambda)=A C_{aging}'(\lambda)$ and
$C_{aging}(\cdot)$ monotone increasing on $[0,1]$
from $C_{aging}(0)=0$ to $C_{aging}(1)=q(\b)$, proposed in \cite{CK}
for the case in which $q(\b_c)>0$, thus allows one to 
set the positive constant $A=\IJ_{\gamma(\b)}/(\b^2 q(\b) \nu'(q(\b))$
in terms of the FDT solution. Unfortunately, even under this 
ansatz, the form of $C_{aging}(\cdot)$ is yet unclear.

Finally, note that it is crucial for our analysis to have 
$\nu'(0)=0$, an assumption we make throughout this work.
That is, our analysis holds in the absence of a random magnetic field.

\section{Limiting exactly spherical dynamics}


Proposition 1.3 and (2.13) of \cite{BDG2} show that 
$\sup_{t \geq 0} K(t) < \infty$. Further, as $C(s,t)$ 
is the limit of the empirical correlation functions $C_N (s,t)$
it is a non-negative definite kernel
on $\reals_+ \times \reals_+$ and in particular, $C(s,t)^2 \leq K(s) K(t)$,
whereas since $\chi(s,t)=\int_0^t R(s,u) du$ is the limit of 
$\chi_N(s,t)$, it follows from the definition of $\chi_N(s,t)$ that 
\begin{equation}\label{eq:rbd}
|\int_{t_1}^{t_2} R(s,u) du |^2 \le K(s) (t_2-t_1) \,,   
\qquad 0 \leq t_1 \leq t_2  \leq s < \infty \,.
\end{equation}

To complete the proof of 
Proposition \ref{prop-sphere}, 
we first prove that
any solution $(R,C,K)$ of \eqref{eqR}--\eqref{eqZ}
consists of positive functions, a key fact in our forthcoming
analysis.
\begin{lem}\label{pos}
For any $f:\R_+ \ra\R$ whose derivative is bounded above
on compact intervals and any $K(0)>0$,  a solution $(R,C,K)$
to \eqref{eqR}--\eqref{eqZ}, if it exists, is positive at all times.
\end{lem}
\prf
By definition $K(t)\ge 0$ for all $t\in\R_+$. Suppose that 
$$
S=\inf\{ u \geq 0: \, C(t,u)\le 0 \mbox{ for some } t \leq u \} < \infty\,.
$$
By continuity of $(C,K)$, since $K(0)>0$ also $S>0$.
Set $\L(s,t)=\exp(-\int_t^s \mu(u) du) >0$ for
$\mu(u)=f'(K(u))$ which is bounded above on compact intervals,
and $R(s,t)=\L(s,t) H(s,t)$. Then, by \cite{GM}, for $s\ge t$,
\begin{equation}
\label{gmfor}
H(s,t) = 1+\sum_{n\ge 1}\beta^{2n}\sum_{\s\in \NC_n}
\int_{t\le t_1\cdots \le t_{2n}\le s}
\prod_{i\in \cro(\sigma) } \nu''(C(t_i,t_{\s(i)})) \prod_{j=1}^{2n} dt_j
\end{equation}
where  $\NC_n$ denotes the
set of involutions of $\{1,\cdots,2n\}$
without fixed points and without crossings
and   $\cro(\sigma)$
is defined to be the set of indices $1\le i\le 2n$
such that $i<\sigma(i)$.
Consequently, 
$$
R(s,t)\ge \L(s,t) > 0 \mbox{ for }t\le s\le S\,,$$
and thus, (\ref{eqC}) implies that also
$$
C(s,t)\ge K(t) \L(s,t)
\mbox{ for }t\le s\le S \,.
$$
Note that in the last two estimates we used the fact that
$\nu'(\cdot)$ and $\nu''(\cdot)$ are non negative on $\R_+$.
Similarly, from the equation \eqref{eqZ} we see that
$\partial_s [\L(s,0)^{-2} K(s)] \ge \L(s,0)^{-2}$ for all $s\le S$
resulting with
$$
K(s)\ge K(0) \L(s,0)^2 + \int_0^s \L(s,v)^2 dv  >0
$$
Hence, the continuous functions $R(s,t),C(s,t)$ and $K(s)$ are
bounded below by a strictly positive constant 
for $t\le s\le S$ in contradiction with the definition of $S$.
We thus deduce that $S=\infty$ and by the preceding argument
both $R(s,t)$ and $C(s,t)$ are positive functions.
\hfill\qed

We next show that if $(R_L,C_L, K_L)$
are solutions of the system (\ref{eqR})--(\ref{eqZ})
with potential $f_L(\cdot)$ as in (\ref{eq:fdef}),
then $K_L(s) \to 1$ as $L \to \infty$, uniformly
over compact intervals. Specifically,
\begin{lem}\label{lem-bdd}
Assuming $K_L(0)=1$, we have that $K_L(s)\ge 1$ 
for all $L>0$ and $s \geq 0$.
Further, for any $T$ finite there exists $B(T)<\infty$, such that 
$K_L(s)\leq 1+B(T) L^{-1}$ for all $s\le T$ and $L \geq B(T)$.
\end{lem}
\prf We first deal with the lower bound on $K_L(\cdot)$.
To this end, fix $L>0$ and let $g(x):=1- 2x f_L'(x)$. 
It is easy to check that $g(x) = 1 + 4L x(1-x)- x^{2k} \ge 0$
for $x\in [0,1]$ and $g(x)\le 0$ for $x \geq 1$.
By Lemma \ref{pos}, we know that the functions $R_L(\cdot,\cdot)$,
and $C_L(\cdot,\cdot)$ are non negative, as is $\psi(r)$ for $r \geq 0$,
so from (\ref{eqZ}) we get the lower bound $\partial_s K_L(s)\ge g(K_L(s))$.
Thus, with $\phi(x)$ a differentiable function 
that is strictly increasing on $[0,1]$ and such that $\phi(x)=1$ for 
all $x \geq 1$, we find that
$$
\partial_s \phi(K_L(s))\ge \phi'(K_L(s)) g(K_L(s))\ge 0 \,.
$$
Consequently, $\phi(K_L(s)) \ge \phi(K_L(0))=\phi(1)=1$ for all $s \geq 0$,
implying by the choice of $\phi(\cdot)$ 
that $K_L(s)\ge 1$ for all $s\ge 0$.

Turning now to the complementary upper bound,
recall that $\psi(r)$ is a polynomial of degree $m-1$, hence
there exists $\kappa<\infty$ such that
$\psi(a b) \leq \kappa (1+a^2)^{m/2} (1+b^2)^{m/2}$ for all $a,b$.
Thus, by (\ref{eq:rbd}), the monotonicity of $\psi(r)$ on $\R_+$
and the non-negative definiteness of
$C_L(s,u)$ we have that for any $s,t,u \geq 0$,
$$
\psi(C_L(s,u))\le \kappa (1+K_L(u))^{\frac{m}{2}}
(1+K_L(s))^{\frac{m}{2}},\qquad \int_0^t R_L(s,u)
du\le \sqrt{t K_L(s)},
$$
and from (\ref{eqZ}) we find that
\begin{equation}\label{eq:Klbd}
\partial_s K_L(s) \leq g(K_L(s)) + 2 \beta^2 \kappa
(1+\sup_{u \leq s} K_L(u))^m \sqrt{K_L(s)} \sqrt{s}  \,.
\end{equation}
Setting now $B(T)=1+ 4 \beta^2 \kappa 3^m \sqrt{T}$ and
fixing $T<\infty$ and $L \geq B(T)$, let
$$
\tau := \inf\{ u\ge 0: K_L (u)\ge 1+B L^{-1}\} \,.
$$
By the continuity of $K_L(\cdot)$ and the fact that
$K_L(0)=1<1+B L^{-1}$, we have that 
$\tau>0$ and further, if $\tau<\infty$ then necessarily
$$
K_L(\tau)=\sup_{u \leq \tau} K_L(u) = 1+B L^{-1} \leq 2 \,.
$$
Recall that $g(x) \leq 1 + 4L(1-x)$ when $x \geq 1$, whereas
from (\ref{eq:Klbd}) we see that if $\tau<\infty$ then
$$
\partial_s K_L (s)_{\mid_{s=\tau}} \le 1 - 4 B + 
4 \beta^2 \kappa  3^m \sqrt{\tau} \,.
$$
Recall the definition of $\tau<\infty$ implying that
$\partial_s K_L (s) \ge 0$ at $s=\tau$. Hence,
our choice of $B=B(T)$ results with $\tau>T$. That is,
$K_L(s) \leq 1+B L^{-1}$ for all $s \leq T$ and $L \geq B(T)$,
as claimed. 
\hfill\qed

Let $\mu_L(s)=f'_L(K_L(s))$ and $h_L(s)=\partial_s K_L(s)$.
Fixing hereafter $T<\infty$, we next prove the equicontinuity and 
uniform boundedness of $(C_L,R_L,K_L,\mu_L,h_L)$, en-route to 
having limit points for $(C_L,R_L,K_L)$.
\begin{lem}\label{lem-tight}
The continuous functions $(C_L(s,t),R_L(s,t))$ and their derivatives 
are bounded uniformly in $L \geq B(T)$ and $0 \leq t \leq s \leq T$. 
Further, the continuous functions $(\mu_L(s),h_L(s))$ and their 
derivatives are bounded uniformly in $L \geq B(T)$ and $s \in [0,T]$.
\end{lem}
\prf Recall that by Lemma \ref{lem-bdd}, for any $L \geq B(T)$,
\begin{equation}\label{eq:kbdd}
\sup_{s\le T}|K_L(s)-1|\le \frac{B(T)}{L} \,.
\end{equation}
Consequently, the collection 
$\{C_L(s,t), 0 \leq t \leq s \leq T, L\geq B\}$ is uniformly bounded.
By (\ref{eq:kbdd}) and our choice of $f_L(r)$, we have that
$$
|\mu_L(s)| \leq 2L |K_L(s)-1| + K_L(s)^{2k-1} \leq 2 B(T) + 2^{2k-1}\,,
\qquad \forall L \geq B(T), s \leq T \,,
$$
whereas by \eqref{gmfor} the collection 
$\{R_L(s,t), 0 \leq t \leq s \leq T, L \geq B\}$ is also uniformly bounded. 
Further, since
\begin{equation}\label{hlid}
h_L(s)=1-2 K_L(s) \mu_L(s) + 2\beta^2 \int_0^s \psi(C_L(s,u)) R_L(s,u) du \,,
\end{equation}
it follows from the uniform boundedness of
$K_L$, $\mu_L$, $C_L$ and $R_L$ 
that $\{h_L(s), s \in [0,T], L \geq B\}$ is also uniformly 
bounded. By the same reasoning, from \eqref{eqR} and \eqref{eqC} 
we deduce that $\partial_s C_L(s,t)$ and $\partial_s R_L(s,t)$
are bounded uniformly in $L \geq B$ and $s,t \in [0,T]$.

Next, differentiating the identity \eqref{gmfor} with respect to $t$, 
we get for $f=f_L$ that
$$
\partial_t H_L (s,t) =
\sum_{n\ge 1}\beta^{2n}\sum_{\s\in \NC_n}
\int_{t=t_1 \le t_2\cdots \le t_{2n}\le s}
\prod_{i\in \cro(\sigma) } \nu''(C_L(t_i,t_{\s(i)})) \prod_{j=2}^{2n} dt_j \,,
$$
where $NC_n$ denotes the finite set of non-crossing involutions 
of $\{1,\ldots,2n\}$ without fixed points.
With the Catalan number $|NC_n|$ bounded by $4^n$, and
since $C_L(t_i,t_{\s(i)}) \in [0,2]$ 
for $t_i, t_{\s(i)}\le T$, $L \geq B(T)$, we thus deduce 
by the monotonicity of $r \mapsto \nu''(r)$ that 
$$
0 \leq \partial_t H_L(s,t) \leq \sum_{n\ge 1}\frac{\beta^{2n}}{(2n-1)!}
4^n \nu''(2)^n (s-t)^{2n-1} \,,
$$
so $\partial_t H_L(s,t)$ is finite and 
bounded uniformly in $L \geq B(T)$ and $0 \leq t \leq s \leq T$. Since
$$
\partial_t R_L (s,t) = \mu_L(t) R_L(s,t) - 
e^{-\int_t^s \mu_L(u) du} \partial_t H_L (s,t) \,,
$$
we thus have that $|\partial_t R_L (s,t)|$ is also bounded uniformly 
in $L \geq B(T)$ and $0 \leq t \leq s \leq T$.

In the course of proving \cite[Lemma 4.1]{BDG2}, one finds
that $\partial_t C(s,t) = R(s,t)+D(s,t)$ for the function $D(s,t)$ of 
\cite[(4.2)]{BDG2}. Consequently, in case of $f=f_L$ we have that 
$$
\partial_t C_L(s,t) = R_L(s,t) - \mu_L(t) C_L (s,t) +
\b^2 \int_0^t C_L (s,u) R_L (t,u) \nu''(C_L(t,u)) du
+ \b^2 \int_0^s \nu'(C_L(t,u)) R_L(s,u) du \,,
$$
which of course is also 
bounded uniformly in $L \geq B(T)$ and $0 \leq t \leq s \leq T$.

Turning to deal with $h_L(\cdot)$, setting   
$g_L(r):=[f_L'(r) r]'-2L = 4L(r-1)+k r^{2k-1}$, 
we deduce from \eqref{eq:kbdd} that 
$|g_L(K_L(s))| \leq 4 B(T) + k 2^{2k}$ for any $s \leq T$ 
and $L \geq B(T)$. Differentiating \eqref{hlid} we find that
$\partial_s h_L(s)= - 4 L h_L(s) +\kappa_L(s)$ for 
$$
\kappa_L(s)= - 2 g_L (K_L(s)) h_L(s)
+ 2\beta^2 \partial_s \int_0^s \psi(C_L(s,u)) R_L(s,u) du \,,
$$
which is thus bounded uniformly in $L \geq B(T)$ and $s \leq T$ 
(in view of the uniform boundedness of $h_L$, $C_L$, $R_L$, 
$\partial_s C_L$ and $\partial_s R_L$). Further, 
recall that $K_L(0)=1$, so by \eqref{eqZ} and our choice of $f_L(r)$
we have that $h_L(0)=1-2f_L'(1)=0$, resulting with 
$$
h_L(s) = \int_0^s e^{-4 L (s-u)} \kappa_L(u) du \;.
$$ 
Since $A(T) = \sup \{ |\kappa_L(u)| : L \geq B(T), u \leq T\}$ is finite, 
we thus deduce that for $L \geq B(T)$,
\begin{equation}\label{eq:hbd}
\sup_{s \leq T} |h_L(s)|\le \frac{A(T)}{4 L} \,,
\end{equation}
and the uniform boundedness of $|\partial_s h_L(s)|$ follows.

Finally, by definition, 
$\partial_s \mu_L(s) = f''_L(K_L(s)) h_L(s)$, yielding for our choice
of $f_L$ that
$$
|\partial_s \mu_L(s)| \leq (2L+(2k-1) 2^{2k-3}) |h_L(s)| \,,
\qquad \forall L \geq B(T), s \leq T \,,
$$
which by (\ref{eq:hbd}) provides 
the uniform boundedness of $|\partial_s \mu_L(s)|$.     
\hfill\qed 

\medskip
\noindent{\bf Proof of Proposition \ref{prop-sphere}.}
By Lemma \ref{lem-tight} we have that  
$(C_L(s,t),R_L(s,t))$, $L \geq B(T)$ are equicontinuous and 
uniformly bounded on $0 \leq t \leq s \leq T$. Further, 
$(K_L,\mu_L,h_L)$ are then equicontinuous and uniformly bounded
on $[0,T]$. By the Arzela-Ascoli theorem, the collection 
$(C_L,R_L,K_L,\mu_L,h_L)$ thus has a limit point 
$(C,R,K,\mu,h)$ with respect to uniform convergence on
$[0,T]$ (or $0 \leq t \leq s \leq T$, whichever is relevant).

By Lemma \ref{lem-bdd} we know that the limit $K(s)=1$
for all $s \leq T$, whereas by \eqref{eq:hbd} we have that
$h(s)=0$ for all $s \leq T$. Consequently, considering 
\eqref{hlid} for the subsequence $L_n \to \infty$ for which 
$(C_{L_n},R_{L_n},K_{L_n},\mu_{L_n},h_{L_n})$ converges to 
$(C,R,K,\mu,h)$ we find that the latter must satisfy 
\eqref{eqZs}. Further, since $R_L(t,t)=1$ and
$C_L(t,t)=K_L(t)$, integrating \eqref{eqR} and \eqref{eqC}
we find that 
$R_L(s,t)= 1 + \int_t^s A_L(\theta,t) d\theta$
and $C_L(s,t)=K_L(t)+\int_t^s B_L(\theta,t) d\theta$, for
\begin{eqnarray*}
A_L(\theta,t) &=& - \mu_L(\theta) R_L(\theta,t) + \b^2 \int_t^\theta
 R_L(u,t) R_L(\theta,u) \nu''(C_L(\theta,u)) du , \\
B_L(\theta,t) &=& - \mu_L(\theta) C_L(\theta,t) + \b^2 \int_0^\theta
  C_L(u,t) R_L(\theta,u) \nu''(C_L(\theta,u)) du
+ \b^2 \int_0^t \nu'(C_L(\theta,u)) R_L(t,u) du \,.
\end{eqnarray*}
Since $A_{L_n}(s,t)$ and $B_{L_n}(s,t)$ converge uniformly
on $0 \leq t \leq s \leq T$ to the right-hand-sides of 
\eqref{eqRs} and \eqref{eqCs}, respectively, we deduce that
for each limit point $(C,R,\mu)$, the functions $C(s,t)$ and $R(s,t)$
are differentiable in $s$ on $0 \leq t \leq s \leq T$ and all 
limit points satisfy the equations \eqref{eqRs}--\eqref{eqZs}.
Further, $C_L(s,t)$ are non-negative functions, and also
symmetric non-negative definite kernels with $C_L(t,t) \to 1$.
Consequently, each of their limit points corresponds 
to a $[0,1]$-valued
symmetric non-negative kernel on $[0,T]^2$. Similarly, since
$R_L(t,t)=1$ and $R_L(s,t)$ satisfy \eqref{eq:rbd},  
the same applies for any limit point $R(s,t)$. The latter
are extended to functions on $[0,T]^2$ by setting 
$R(s,t)=R_L(s,t)=0$ whenever $s<t$. 

Finally, it is not hard to verify 
that the system of equations \eqref{eqRs}--\eqref{eqZs}
with $C(s,t)=C(t,s)$, $C(t,t)=R(t,t)=1$ and $R(s,t)=0$ for $s<t$,
admits at most one bounded solution $(R,C)$ on $[0,T]^2$.
Indeed, considering the difference between the integrated form of 
\eqref{eqRs}--\eqref{eqCs} for two such solutions 
$(C,R)$ and $(\bar C,\bar R)$, 
since $\nu''$ is uniformly Lipschitz on $[0,1]$, 
the functions $\Delta R(s,t)=|R(s,t)-\bar R(s,t)|$
and $\Delta C(s,t)=|C(s,t)-\bar C(s,t)|=\Delta C(t,s)$ are such that
for $0\le t\le s\le T$,
\begin{eqnarray}
\Delta R(s,t)&\le& \kappa_1 
[\int_t^s \Delta R(v,t) dv +\int_t^s h(v) dv ]
\label{tutuR}\\
\Delta C(s,t)&\le& \kappa_1 [
\int_t^s \Delta C(v,t) dv +h(t)+\int_t^s h(v) dv]
\label{tutuC}
\end{eqnarray}
where $h(v) := \int_0^v [\Delta R(v,u) +\Delta C(v,u)] du$ and
$\kappa_1 < \infty$ depends on $T$, $\b$, $\nu(\cdot)$ and
the maximum of $|R|$, $|C|$, $|\bar R|$ and $|\bar C|$ on $[0,T]^2$.
Integrating these inequalities over $t \in [0,s]$, 
since $\Delta R(v,u)=0$ for $u \geq v$ 
and $\Delta C(v,u)=\Delta C(u,v)$, we find that 
\begin{eqnarray*}
\int_0^s \Delta R(s,t) dt &\le& \kappa_2 \int_0^s h(v) dv \,, \\
\int_0^s \Delta C(s,t) dt &\le& \kappa_2 \int_0^s h(v) dv \,,
\end{eqnarray*}
for some finite constant $\kappa_2$ (of the same type of 
dependence as $\kappa_1$).
Summing the latter two inequalities we see that for all $s \in [0,T]$,
$$
0 \leq h(s) \leq 2 \kappa_2 \int_0^s h(v) dv. 
$$
Further, $h(0)=0$, so by Gronwall's lemma $h(s)=0$
for all $s \in [0,T]$.  Plugging this result back into
\eqref{tutuR}-\eqref{tutuC} and observing that 
$\Delta R(t,t)=\Delta C(t,t)=0$, we deduce that 
$\Delta R(s,t)=\Delta C(s,t)=0$
for all $0 \le t\le s \le T$, yielding the stated uniqueness.

In conclusion, when $L \to \infty$ the collection $(C_L,R_L,K_L)$ 
converges towards the unique solution $(C,R,K)$ of 
\eqref{eqRs}--\eqref{eqZs}, as claimed.
\hfill \qed

\section{Exponential decay for small values of $\beta$}\label{expdec}
We consider here the solution $(R,C)$ of 
(\ref{eqRs})--(\ref{eqZs}) and prove Proposition \ref{prop-exp-tails}
about its 
exponential decay in $s-t$ for all $\beta>0$ sufficiently small. 

\subsection{Exponential decay of $R(s,t)$ for $\beta$ small}

\begin{lem}\label{expdecR}
If $\beta < 1/(4\sqrt{\nu''(1)})$ then $R(s,t)\le e^{-\d_\b (s-t)}$ for 
$\d_\beta = \frac{1}{2} - 2\beta \sqrt{\nu''(1)}>0$ and all $s\ge t$.
\end{lem}
\prf By \cite{GM}, for $s \geq t$, the solution $(R,C)$ of 
\eqref{eqRs}--\eqref{eqCs} is such that $R(s,t)=\L(s,t) H(s,t)$ for
$\L(s,t)=\exp(-\int_t^s \mu(u) du)$ with
$\mu(u)$ of \eqref{eqZs} and $H(s,t)$ of \eqref{gmfor}.
Further, since $C(s,t) \in [0,1]$ for
all $s,t \geq 0$, it follows that $\nu''(C(s,t)) \leq \nu''(1)<\infty$.
Recall that $|NC_n| \le 4^n$, hence we deduce that
$$
H(s,t)\le 1+\sum_{n\ge 1}\b^{2n} 4^n \nu''(1)^n
\frac{(s-t)^{2n}}{2n!} \le \sum_{k\ge 0}(2\b \sqrt{\nu''(1)})^{k}
\frac{(s-t)^k}{k!}= e^{2 \beta \sqrt{\nu''(1)}(s-t)} \,.
$$
By Lemma \ref{pos} we know that $C$ and $R$ are non-negative functions, 
hence $\mu(u)\ge \frac{1}{2}$ (by \eqref{eqZs}), resulting with
$\L(s,t) \le e^{-(s-t)/2}$. In conclusion,
$$
R(s,t) = \L(s,t) H(s,t) \le e^{(2\beta \sqrt{\nu''(1)}-\frac{1}{2})(s-t)}
= e^{-\delta_\beta (s-t)} \,,
$$
where $\delta_\beta>0$ for 
$\beta < 1/(4\sqrt{\nu''(1)})$, as claimed.
\hfill\qed

\subsection{Exponential decay of $C(s,t)$ for $\beta$ small}
\begin{lem}\label{expdecC} 
For some $\b_0>0$ and any $\beta<\beta_0$ there exist $M_\beta<\infty$
and $0<\eta_\b<\delta_\b$ such that 
$$
C(s,t)\le M_\beta e^{-\eta_\b |s-t|} \,.
$$
\end{lem}
\prf From equation (\ref{eqC}) we get that for $s \geq t \geq 0$,
\begin{equation}\label{eqC2}
C(s,t)= \L (s,t) +\beta^2 \int_t^s \L(s,v) I_1(v,t) dv +
\beta^2 \int_t^s \L(s,v) I_2(v,t) dv
\end{equation}
with $\L(s,v)=\exp(-\int_v^s \mu(u) du)$,
\begin{eqnarray}
I_1(v,t)&=&\int_0^v C(u,t) R(v,u) \nu''(C(v,u)) du \,, \label{eqI1} \\
I_2(v,t)&=& \int_0^t \nu'(C(v,u)) R(t,u) du \,.  \label{eqI2}
\end{eqnarray}
For $\beta < 1/(4\sqrt{\nu''(1)})$ we know from Lemma \ref{expdecR} 
that $R(s,t)\le e^{-\d_\beta(s-t)}$. With $C(v,u)=C(u,v) \in [0,1]$
and  since $\nu'(r)\le \nu''(1) r$ and $\nu''(r) \leq \nu''(1)$
for $0\le r\le 1$, we get that for $v \geq t$, 
\begin{eqnarray*}
I_1(v,t)&\le & \nu''(1) \int_0^v C(u,t) e^{-\d_\beta (v-u)}
du, 
\\
I_2(v,t)&\le &\nu''(1)  \d_\beta^{-1} \sup_{u\le t} C(u,v).
\end{eqnarray*}
Recall that $\mu(u) \geq 1/2$, so $\L(s,v) \le e^{-(s-v)/2}$.
Hence, with the symmetric function
$\D(t,s):=\sup_{u\le t,v \le s} \, C(u,v)$
we deduce from (\ref{eqC2}) that for $s \ge t\ge 0$,
\begin{eqnarray*}
\D(t,s)&\le&  e^{-\frac{1}{2}(s-t)}
+\beta^2 \nu''(1)  \int_t^s e^{-\frac{1}{2}(s-v)}[
\int_0^v C(u,t) e^{-\d_\beta (v-u)}
du + \d_\beta^{-1}\D(t,v)] dv\\
&\le& e^{-\frac{1}{2}(s-t)}
+\beta^2 \nu''(1)  \int_t^s e^{-\frac{1}{2}(s-v)}
\int_0^t  e^{-\d_\beta (v-u)} du dv\\
&&
+\beta^2 \nu''(1)\int_t^s \D(t,v)[\d_\beta^{-1}
e^{-\frac{1}{2}(s-v)}+\int_t^v e^{-\frac{1}{2}(s-v)-\d_\beta (v-u)}
du] dv
\end{eqnarray*}
It is straightforward to see that for any $\d\in (0, 1/2)$ and $s \geq t$,
\begin{equation}\label{eq:integ}
\int_t^s e^{-\frac{1}{2} (s-v)-\d(v-t)} dv \le 2(1-2\d)^{-1}e^{-\d (s-t)}
\end{equation}
and with $\d_\beta \in (0,1/2)$ we thus obtain for $s \geq t$ the bound
\begin{eqnarray*}
\D(t,s)&\le& M_\beta e^{-{\d_\beta} (s-t)} 
+A_\beta \int_t^s \D(t,v) e^{-{\d_\beta} (s-v)} dv \,,
\end{eqnarray*}
with $M_\beta=1+2\beta^2 \nu''(1)(1-2\d_\beta)^{-1}\d_\beta^{-1}$ and
$A_\beta= 2 \beta^2\nu''(1) \d_\beta^{-1}$. Therefore, fixing $t \geq 0$,
the function $h_t(s)=e^{{\d_\beta}( s-t)} \D(t,s)$ satisfies
$$h_t(s)\le M_\beta+ A_\beta\int_t^s h_t(v)dv, \quad s\ge t, $$
and so by  Gronwall's lemma $h_t(s)\le  M_\beta e^{A_\beta(s-t)}$.
We therefore conclude that for any $s \geq t$, 
$$C(s,t)\le  M_\beta e^{-(\d_\beta-A_\beta)( s-t)} \,,
$$
which proves the lemma since for $\beta \to 0$ 
we have that $\d_\beta \to 1/2$ while $A_\beta \to 0$ (and so 
$\eta_\beta=\d_\beta-A_\beta>0$ for any $\beta>0$ small enough).
\hfill\qed 

\section{Getting the FDT equations}\label{fdtequation}
With $\bD=\{(s,t): 0 \leq t \leq s \} \subset \R_+ \ts \R_+$,
we consider the map $\Psi:(R,C)\ra (\tilde R,\tilde C)$ on  
$$
\Aa^+ =\{ (R,C)\in \CC(\bD) \ts \CC (\R_+\ts\R_+) 
: R(t,t)=C(t,t)=1, \; R(s,t) \ge 0,
C(s,t)=C(t,s)\ge 0 \} \,,
$$
such that for $s \geq t$,
\begin{eqnarray}
\partial_s \tilde R(s,t)&\!\!\!\! =\!\!\!\!
& - \mu_{R,C}(s)  \tilde R(s,t) + \b^2 \int_t^s
\tilde R(u,t)\tilde R(s,u) \nu''(C(s,u)) du ,\label{eqRP}\\
\partial_s \tilde C(s,t)&\!\!\!\! = \!\!\!\!
& - \mu_{ R, C}(s) \tilde C(s,t) +
\b^2 \int_0^s C(u,t) R(s,u) \nu''(C(s,u)) du
+ \b^2 \int_0^t \nu'(C(s,u))  R(t,u) du,\label{eqCP}
\end{eqnarray}
where $\tilde R(t,t)=\tilde C(t,t)=1$, $\tilde C(t,s)=\tilde C(s,t)$ and
\begin{equation}\label{eqmuP}
\mu_{ R, C}(s)=\frac{1}{2} +\beta^2\int_0^s \psi( C(s,u)) R(s,u) du.
\end{equation}
Assuming $(R,C) \in \Aa^+$, we have that
$\mu_{R,C}(s) \geq 1/2$ is continuous and further, by \cite{GM} 
there exists a unique non-negative solution $\tilde R(s,t)$
of \eqref{eqRP} which is continuous on $\bD$ (see for example
\eqref{gmfor} for existence, uniqueness and non-negativity of the
solution, and the proof of Lemma \ref{lem-tight} 
for the differentiability, hence 
continuity of $\tilde R(s,t)$). With $C \geq 0$ and $R \geq 0$, clearly there
is also a unique non-negative solution $\tilde C(s,t)$ to
\eqref{eqCP} which is continuous on $\bD$ and due to the boundary
condition $C(t,t)=1$, its symmetric extension to $\R_+\ts \R_+$ remains
continuous. Thus, $\Psi(\Aa^+)\subseteq \Aa^+$.

We proceed to show that for small $\beta$
and a suitable choice of the positive constants $\d,r,\rho,c$
the solution $(R,C)$ of \eqref{eqRs}- \eqref{eqCs}
is a fixed point of the mapping $\Psi$ on the space
$$
\Sa(\d,r,\rho,c)=\{(R,C)\in\Aa(\d,r,\rho,c) :
\exists R_{\fdt}(\tau) = \lim_{t \to \infty} R(t+\tau,t) \,,\;
\exists C_{\fdt}(-\tau) = 
C_{\fdt}(\tau) = \lim_{t \to \infty} C(t+\tau,t) \,,\;
\forall \tau \geq 0 \}\,,
$$
where 
$$
\Aa(\d,r,\rho,c)=\{(R,C)\in \Aa^+ : C(s,t)\le c e^{-\d|s-t|},
\qquad
R(s,t)\le \rho (r|s-t|+1)^{-3/2} e^{-\d (s-t)}
\quad\mbox{ for all } s\ge t\}\,.
$$

This of course implies that 
the solution $(R,C)$ of \eqref{eqRs}- \eqref{eqCs} is
such that the functions 
\begin{eqnarray}
R_{\fdt}(\tau) &=& \lim_{t\ra\infty} R(t+\tau,t) \,, \label{Rfdt-ex}\\
C_{\fdt}(\tau) &=& \lim_{t\ra\infty} C(t+\tau,t) \,, \label{Cfdt-ex}
\end{eqnarray}
are well defined for all $\tau \geq 0$. Further, 
for any $(R,C) \in \Sa(\d,r,\rho,c)$ the corresponding 
functions $(R_{\fdt},C_{\fdt})$ are clearly in the set
\begin{eqnarray*}
\Ba(\d,r,\rho,c):=\{(R,C) \in \Ba(\R_+) \ts \Ba(\R) : 
C(0)=R(0)=1, \;\; &&
0\le C(\tau)=C(-\tau) \le ce^{-\d|\tau|}, \;\; \\ &&
0\le R(\tau)\le \rho (r\tau +1)^{-3/2} e^{-\d \tau}\} \,,
\end{eqnarray*}
so in particular, \eqref{Cfdt-ex} holds for all $\tau \in \R$
and $(R_{\fdt},C_{\fdt})$ of \eqref{Rfdt-ex}--\eqref{Cfdt-ex}
are of exponential decay. 

To this end, we start by finding constants 
$\d,\rho,c$ and $r=r(\b)$ for which $\Sa(\d,r,\rho,c)$ is 
closed under the mapping $\Psi$.
\begin{prop}\label{Psi}
There exist finite, positive universal constants $c_1$ and $M_1 \geq 2$,
such that for $c=2$, $\rho=c_1$, $\d=1/6$ and
$r=\beta\sqrt{\nu''(c)} \leq 1/(3 M_1)$, both
\begin{equation}\label{al2}
\Psi(\Aa(\d,r,\rho, c))\subseteq \Aa(\d,r,\rho, c)\,,
\end{equation}
and 
\begin{equation}\label{al3}
\Psi(\Sa(\d,r,\rho, c))\subseteq \Sa(\d,r,\rho, c).
\end{equation}
\end{prop}

\nn {\bf Proof of Proposition \ref{Psi}:} We first verify that 
\eqref{al2} holds. To this end, 
setting $\tilde R(s,t)= \L(s,t) \tilde H(s,t)$ for
\begin{equation}\label{eqL}
\L(s,t)=e^{-\int_t^s \mu_{R,C}(u) du}\,,
\end{equation}
we have that $\tilde H(t,t)=1$. Further, from \cite{GM} we have that
for any $(s,t) \in \bD$,
\begin{equation}\label{eqtH}
\tilde H(s,t)=1+\sum_{n\ge 1}\beta^{2n} \sum_{\sigma\in\NC_n}
\int_{t\le t_1\cdots\le t_{2n}\le s}
\prod_{i\in \cro(\sigma) } 
\nu''(C(t_i,t_{\sigma_i}))\prod_{j=1}^{2n} dt_j \,.
\end{equation}
Consequently, as $|\NC_n|= (2\pi)^{-1}\int_{-2}^2 x^{2n} \sqrt{4-x^2} dx$ 
and $C(u,v) \in [0,c]$ for $(R,C)\in \Aa(\d,r,\rho,c)$, we have the bound
\begin{eqnarray}\label{eqtHbd1}
\tilde H(s,t)&\le& \sum_{n\ge 0} (\beta^2 \nu''(c))^n
\sum_{\s\in\NC_n}
\int_{t\le t_1\le\cdots\le t_{2n}\le s}  \prod_{j=1}^{2n} dt_j\\
&=& \sum_{n\ge 0} \frac{(\beta^2 \nu''(c))^n (s-t)^{2n}}{(2n !)}
(2\pi)^{-1}\int_{-2}^{2} x^{2n} \sqrt{4-x^2} dx \nonumber \\
&=& (2\pi)^{-1}\int_{-2}^2 e^{\beta \sqrt{\nu''(c)} (s-t) x} \sqrt{4-x^2} dx \,.
\nonumber
\end{eqnarray}
It is well known (see for example \cite[(3.8)]{2001}) that
for some universal constant $1 \leq c_1 <\infty$ and all $\theta$,
$$
(2\pi)^{-1}\int_{-2}^2 e^{\theta x} \sqrt{4-x^2} dx\le 
c_1(1+|\theta|)^{-3/2}\, e^{2 | \theta |}\,,
$$
from which we thus deduce that 
\begin{equation}\label{eq:tHbd} 
\tilde H(s,t)\le c_1
(1+ \beta \sqrt{\nu''(c)} (s-t))^{-3/2}\, e^{2\beta \sqrt{\nu''(c)}(s-t)}\,.
\end{equation}
Further, since $(R,C) \in \Aa^+$, we know that $\mu_{R,C}(u)\ge 1/2$ 
resulting with $\L(s,t)\le e^{-(s-t)/2}$. It then follows that for
our choice of $\rho=c_1$, $r=\beta\sqrt{\nu''(c)}$ and $\d=1/6 \leq 1/2-2r$,
\begin{eqnarray*}
\tilde R(s,t) &\le& c_1(1+ \beta\sqrt{\nu''(c)} (s-t) )^{-3/2}\,
e^{-(\frac{1}{2}-2 \beta\sqrt{\nu''(c)}) (s-t)} \\
&\leq& \rho (1+r (s-t))^{-3/2} e^{-\d (s-t)} \,.
\end{eqnarray*}

Considering next the function 
$\tilde C$, recall that $\tilde C(t,t)=1$ 
and for $(R,C)\in \Aa(\d,r,\rho,c)$ we have
that $\nu'(C(v,u))\le \nu''(c) C(v,u)$ and 
$\L(s,v)\le e^{-(s-v)/2}$ for all $v \leq s$ and $u$. Thus,
we get from \eqref{eqCP} that for $(s,t) \in \bD$ and
$(R,C)\in \Aa(\d,r,\rho,c)$,
\begin{eqnarray*}
&& \tilde C(s,t)= \L(s,t)+ \b^2 \int_t^s \L(s,v) dv[\int_0^v  C(u,t)
 R(v,u) \nu''(C(v,u)) du +  \int_0^t  \nu'(C(v,u)) R(t,u) du]\\
&\le& 
\L(s,t) + \b^2 c \rho \nu''(c) \int_t^s e^{-\frac{1}{2}(s-v)}
e^{-\d( v-t)}  dv
[\int_0^v (r(v-u)+1)^{-\frac{3}{2}} du
+\int_0^t (r(t-u)+1)^{-\frac{3}{2}} du] \\
&\le& [ 1+ 2 K r^{-1} \b^2 c \rho \nu''(c) \L(s,t) (\frac{1}{2} -\d)^{-1} ]
e^{-\d(s-t)}
\end{eqnarray*}
where in the last inequality we have used \eqref{eq:integ} for $\d=1/6<1/2$
and
$$
K:= \int_0^\infty (\theta+1)^{-\frac32} d\theta = 2 \,.
$$
This shows that $(\tilde R, \tilde C) \in \Aa(\d,r,\rho,c)$ since our
choices of $c=2$, $\d=1/6$, $\rho=c_1$ and $r = \beta \sqrt{\nu''(c)}$,
result with
$$
1+ \frac{2 K c \rho \b^2 \nu''(c)}{(\frac{1}{2} -\d)r}
= 1+ 12 K c_1 r \le 2 = c \,,
$$
once we take $M_1=4 K c_1$ and $\beta$
sufficiently small for $r=r(\b) \leq 1/(3 M_1)$.

\bigskip
Our next task is to verify that \eqref{al3} holds. That is,
assuming that $(R,C) \in \Sa(\d,r,\rho,c)$ we are to show that the 
limits $(\tilde R_{\fdt},\tilde C_{\fdt})$ exist for the solution 
$(\tilde R,\tilde C)$ of \eqref{eqRP}--\eqref{eqmuP}. To this end,
recall that by \eqref{eqL}, \eqref{eqtH}, and 
\eqref{eqC2}, for any $t \geq 0$ and $\tau \geq v \geq 0$, 
\begin{eqnarray*}
\L(t+\tau,t+v) &=& e^{-(\tau-v)/2} e^{-\b^2 \int_v^\tau I_0(t+u,t) du} \\
\tilde R(t+\tau,t)&=& \L(t+\tau,t) \tilde H(t+\tau,t) \\
\tilde H(t+\tau,t)&=& 1 + 
\sum_{n\ge 1}\beta^{2n} \sum_{\s\in\NC_n}
\int_{0\le \theta_1\le \cdots\le \theta_{2n}\le \tau}
\prod_{i\in \cro(\s)} \nu''(C(t+\theta_i,t+\theta_{\sigma(i)}))
\prod_{j=1}^{2n} d\theta_j \\
\tilde C(t+\tau,t)&=& \L(t+\tau,t) +
\b^2 \int_0^{\tau} \L (t+\tau,t+v) I_1(t+v,t) dv
+ \b^2 \int_0^\tau \L (t+\tau,t+v) I_2(t+v,t) dv
\end{eqnarray*}
where by \eqref{eqmuP}, \eqref{eqI1} and \eqref{eqI2},
\begin{eqnarray*}
 I_0(t+v,t)&=& \int_{-t}^v \psi(C(t+v,t+u)) R(t+v,t+u) du \\
 I_1(t+v,t)&=& \int_{-t}^v C(t+u,t) R(t+v,t+u) \nu''(C(t+v,t+u)) du \\
 I_2(t+v,t)&=& \int_{-t}^0 \nu'(C(t+v,t+u)) R(t,t+u) du \,. 
\end{eqnarray*}
Since $\psi(\cdot)$, $\nu''(\cdot)$ and $\nu'(\cdot)$ are 
continuous and $(R,C) \in \Sa(\d,r,\rho,c)$, 
as $t \to \infty$ the bounded integrands in the preceding formulas 
converge pointwise (per fixed $u=v-\theta$) to the corresponding 
expression for $(R_{\fdt},C_{\fdt})$.
Further, by the exponential tails of $(R,C) \in \Aa(\d,r,\rho,c)$ 
the integrals over $[-t,-m]$ in the formulas for $I_0$, $I_1$ and $I_2$,
are uniformly in $t$ bounded by $\rho \d^{-1} \psi(c) e^{-\d m}$.
Thus, applying bounded convergence theorem for the integrals 
over $[-m,v]$, then taking $m \to \infty$, we deduce that
for each fixed $v \geq 0$,
\begin{eqnarray}
\hh{I}_0 &:=& \lim_{t\ra\infty} I_0(t+v,t) =
\int_0^\infty \psi(C_{\fdt}(\theta)) R_{\fdt}(\theta) d\theta, \nonumber
\\
\hh{I}_1(v) &:=&
\lim_{t\ra\infty} I_1(t+v,t) 
= \int_0^\infty
C_{\fdt} (v-\theta) R_{\fdt} (\theta) \nu''(C_{\fdt} (\theta)) d\theta,
\label{eqI1l}
\\
\hh{I}_2(v) &:=&
\lim_{t\ra\infty} I_2(t+v,t) 
=\int_v^\infty \nu'(C_{\fdt} (\theta)) R_{\fdt} (\theta-v)d\theta.
\label{eqI2l}
\end{eqnarray}
By the preceding discussion we also know that
$0 \leq I_i(t+v,t) \leq \rho \psi(c) \d^{-1}$ for
$i=0,1,2$ and all $v,t \geq 0$. Thus, with 
$\L(t+\tau,t+v) \in [0,1]$, by bounded convergence 
for each $\tau \geq v \geq 0$, 
\begin{eqnarray*}
\hh{\L}(\tau-v) &:=& \lim_{t \to \infty} 
\L(t+\tau,t+v) = e^{-(\tau-v)(1/2+\b^2 \hh{I}_0)}\,, \\
\tilde C_{\fdt} (\tau) &:=& \lim_{t \to \infty} 
\tilde C(t+\tau,t)= \hh{\L}(\tau) +
\b^2 \int_0^{\tau} \hh{\L}(\tau-v) \hh{I}_1(v) dv
+ \b^2 \int_0^\tau \hh{\L}(\tau-v) \hh{I}_2(v) dv \,.
\end{eqnarray*}
We also have that for any $n\in\N$, all $\s\in\NC_n$ and 
each fixed $\theta_1,\ldots,\theta_{2n} \geq 0$,
$$
\lim_{t\ra\infty}
\prod_{i\in\cro(\sigma)}\nu''(C (t+\theta_i, t+\theta_{\sigma(i)})) =
\prod_{i \in\cro(\sigma)}\nu''(C_{\fdt}(\theta_i- \theta_{\sigma(i)}))\,,
$$
By bounded convergence, we have the convergence of the 
corresponding integrals over
$0 \leq \theta_1 \leq \cdots \leq \theta_{2n} \leq \tau$.
Further, the  
non-negative series \eqref{eqtH} is dominated in $t$ by a summable 
series (see \eqref{eqtHbd1}), so by dominated convergence,
$$
\tilde H_{\fdt} (\tau) := 
\lim_{t \to \infty} \tilde H(t+\tau,t) = 
1 + \sum_{n\ge 1}\beta^{2n} \sum_{\s\in\NC_n}
\int_{0\le \theta_1\le \cdots\le \theta_{2n}\le \tau}
\prod_{i\in \cro(\s)} \nu''(C_{\fdt}(\theta_i-\theta_{\sigma(i)}))
\prod_{j=1}^{2n} d\theta_j \,.
$$
It thus follows that
$$
\tilde R_{\fdt} (\tau) := \lim_{t \to \infty} 
\tilde R(t+\tau,t)= \hh{\L}(\tau) \tilde H_{\fdt} (\tau) \,,
$$
exists for each $\tau \geq 0$, which establishes our claim \eqref{al3}
(we have already shown that $\tilde C_{\fdt} (\tau)$ exists).
\hfill\qed

\bigskip

We next show that $\Psi$ is a contraction on 
$\Sa(\d,r,\rho,c)$ and provide the set of equations that characterizes 
the functions $R_{\fdt}$ and $C_{\fdt}$ of Theorem \ref{FDT}.
\begin{prop}\label{FDT1}
For $\d,\rho,c$ and $r=r(\beta)$ of Proposition \ref{Psi},
if $\beta$ is small enough then
$\Psi$ is a contraction on $\Aa(\d,r,\rho,c)$ equipped with the norm
\begin{equation}\label{eq:unorm}
\| (R,C)\|= \sup_{(s,t) \in \bD} |R(s,t)|+ \sup_{s,t \in \R_+} |C(s,t)| \,,
\end{equation}
and the solution $(R,C)$ of \eqref{eqRs}--\eqref{eqCs} is also 
the unique fixed point of $\Psi$ in $\Sa(\d,r,\rho,c)$. Consequently, 
the functions $(R_{\fdt},C_{\fdt})$ of \eqref{Rfdt-ex}--\eqref{Cfdt-ex}
are then the unique solution in $\Ba(\d,r,\rho,c)$ of the FDT equations  
\begin{eqnarray}
R'(\tau)&\!\!\!\! =\!\!\!\!
& - \mu R(\tau) + \b^2 \int_0^\tau R(\tau-\theta) 
R(\theta) \nu''(C(\theta))d\theta ,\label{eqRFDT}\\
C'(\tau)&\!\!\!\! = \!\!\!\!
& - \mu C(\tau) 
+\b^2 \int_0^\infty C (\tau-\theta) R (\theta) \nu''(C (\theta)) d\theta
+ \b^2 \int_\tau^\infty \nu'(C(\theta)) R(\theta-\tau)d\theta + \IJ,
\label{eqCFDT}\\
\mu&= & \frac 12 + \b^2
\int_0^\infty \psi(C(\theta)) R(\theta) d\theta + \IJ \,, \label{eqKFDT}
\end{eqnarray}
for $\IJ=0$.
\end{prop}

\noindent {\bf Proof of Proposition \ref{FDT1}:}
Keeping $\d$, $\rho$, $c$ and $r=r(\b)$ as in Proposition \ref{Psi},
we first check that for any $\beta$ small enough, 
$\Psi$ is a contraction on $\Aa(\d,r,\rho,c)$ equipped with the
uniform norm $\|((R,C)\|$ of \eqref{eq:unorm}. To this end, we
consider the pairs $(\tilde R_i,\tilde C_i)=\Psi(R_i,C_i)$ for
$(R_i,C_i) \in \Aa(\d,r,\rho,c)$, $i=1,2$. 
Denoting hereafter in short $\D f(s,t)=|f_1(s,t)- f_2(s,t)|$  
and $\D f(s)=\sup_{u\le v\le s} \D f(v,u)$ when $f$ is
one of the functions of interest to us, such as $R$, $C$,
$\L$, $\tilde H$, $\tilde R$ or $\tilde C$, we shall show that
there exist finite constants $M_R=M_R(\d,\rho,c)$ and $M_C=M_C(\d,\rho,c)$ 
such that for any finite $s$,
\begin{eqnarray}
\D \tilde R (s) &\leq& M_R \b^2 [ \D R(s) + \D C(s) ] \,, 
\label{eq:lipR}
\\
\D \tilde C (s) &\leq& M_C \b^2 [ \D R(s) + \D C(s) ] \,.
\label{eq:lipC}
\end{eqnarray}
So, if $\b$ is small enough for
$M_R \b^2 \leq 1/3$ and $M_C \b^2 \leq 1/3$,
then from \eqref{eq:lipR} and \eqref{eq:lipC} we deduce that
$$
\| (\tilde R_1,\tilde C_1) - (\tilde R_2,\tilde C_2) \| 
= \sup_{s \geq 0} \D \tilde R(s) + \sup_{s \geq 0} \D \tilde C(s)
\leq \frac{2}{3} [\sup_{s \geq 0} \D R(s) + \sup_{s \geq 0} \D C(s)]
= \frac{2}{3} \| (R_1,C_1) - (R_2,C_2) \| \,.
$$
In conclusion, the mapping 
$\Psi$ is then a contraction on $\Aa(\d,r,\rho,c)$, since 
\begin{equation}\label{eq:Psicont}
\| \Psi(R_1,C_1) - \Psi(R_2,C_2) \| 
\leq \frac{2}{3} \| (R_1,C_1) - (R_2,C_2) \| \,,
\end{equation}
for any $(R_i,C_i) \in \Aa(\d,r,\rho,c)$, $i=1,2$.

The challenge in deriving 
\eqref{eq:lipR}--\eqref{eq:lipC} is to get bounds that  
are uniform over $(s,t) \in \bD$. In doing so, we use
the tail estimates for $(R_i(s,t),C_i(s,t))$, $i=1,2$
(which hold for all functions in $\Aa(\d,r,\rho,c)$),
in order to improve upon the arguments of
\cite[proof of Proposition 4.2]{BDG2}, where the Lipschitz
bounds are derived for finite time intervals.

\nn
$\bu$
{\it The Lipschitz bound \eqref{eq:lipR} on $\tilde R$.}

\nn
We rely on the formulas \eqref{eqtH}
and $\tilde R(s,t)=\tilde H(s,t) \L(s,t)$.
Indeed, since $C_1$ and $C_2$ are $[0,c]$-valued
symmetric functions, $t_i \in [0,s]$ and 
both $\nu''(\cdot)$ and $\nu'''(\cdot)$ are 
non-negative and monotone non-decreasing, it follows that
for any $n$, $t_{2n} \leq s$ and $\sigma \in\NC_n$,
$$
| \prod_{i\in \cro(\sigma) } \nu''(C_1(t_i,t_{\sigma_i}))
- 
 \prod_{i\in \cro(\sigma) } \nu''(C_2(t_i,t_{\sigma_i})) |
\leq n \nu''(c)^{n-1} \nu'''(c) \D C(s) \,.
$$ 
Thus, with $r=\beta \sqrt{\nu''(c)}$ 
we easily deduce from \eqref{eqtH} that
\begin{eqnarray}
\D\tilde H(s,t)
&\le&4\beta^{2} \nu'''(c)(s-t)^2
\sum_{n\ge 1}n(2n!)^{-1} [2 r (s-t)]^{2(n-1)}
\D C(s)
\nonumber\\
&\le&  2\beta^2  \nu'''(c) (s-t)^2 e^{2 r (s-t)} \D C(s) \,.
\label{DH}
\end{eqnarray}

Next, note 
that $\psi'(\cdot)$ is a polynomial of non-negative coefficients,
hence non-decreasing on $\R_+$. Consequently, for any
$(R_i,C_i) \in \Aa(\d,r,\rho,c)$ and all $(v,u) \in \bD$,
$$
\D (\psi(C) R) (v,u) \leq \psi'(c) [R_1(v,u) \D C(v,u) + C_2(v,u) \D R(v,u)]
\leq \psi'(c) [\rho \D C(v,u) + c \D R(v,u)] e^{-\d (v-u)} \,.
$$
As $|e^{-x}-e^{-y}| \leq |x-y|$ for all $x,y \geq 0$, we thus get that 
\begin{eqnarray}
\D \L(s,t) &\le & \beta^2 e^{-(s-t)/2} \int_t^s \int_0^v 
 \D (\psi(C) R) (v,u) du dv \nonumber\\
&\le& \beta^2 e^{-(s-t)/2}\psi'(c) [\rho \D C (s) +c \D  R(s)] 
\int_t^s \int_0^v e^{-\d (v-u) } du dv \nonumber\\
&\le & c_2 \beta^2 (s-t) e^{-(s-t)/2} [\D C(s)+\D R(s)] \,,
\label{DL}
\end{eqnarray}
for $c_2=\psi'(c)(\rho+c)\d^{-1}$. Since 
$\tilde R (s,t)=\tilde H(s,t) \L(s,t)$ we now obtain from
\eqref{eq:tHbd}, \eqref{DH} and \eqref{DL} that
\begin{eqnarray*}
\D \tilde R(s,t)&\le& \L_1(s,t) \D \tilde H(s,t)+ \tilde
H_2(s,t) \D \L(s,t) \\
&\le& 2\beta^2  \nu'''(c)  (s-t)^2 e^{-[\frac 12- 2r](s-t)} \D C(s) 
+ c_1 c_2 \beta^2 (s-t) e^{-[\frac 12- 2 r](s-t)} [\D C(s)+\D R(s)]
\end{eqnarray*}
Further, as $1/2-2r \geq 1/6$, the latter bound leads to 
\eqref{eq:lipR} for the finite universal constant 
$$
M_R := \sup_{\theta \geq 0} \,
(2 \nu'''(c) \theta^2 e^{-\theta/6} + c_1 c_2 \theta e^{-\theta/6}) \,.
$$

\smallskip\nn
$\bu${\it The Lipschitz bound \eqref{eq:lipC} on $\tilde C$.}

\nn
Recall that $\tilde C(s,t)$ for $(s,t) \in \bD$ is given by
\eqref{eqC2}, whereas for all 
$(R,C) \in \Aa(\d,r,\rho,c)$ and $0 \leq t \leq v$, 
\begin{eqnarray*}
I_1(v,t) = \int_0^v C(u,t) R(v,u) \nu''(C(v,u)) du &\le& 
c \rho K r^{-1} \nu''(c) e^{-\d(v-t)},\\
I_2(v,t) = \int_0^t \nu'(C(v,u)) R(t,u) du&\le&
 c \rho K r^{-1} \nu''(c) e^{-\d(v-t)} .
\end{eqnarray*}
Further, with $c_3 = \rho(\nu''(c) + c \nu'''(c))$,
$c_4= c \nu''(c)$, $c_5=\d^{-1} \max(c_3,2c_4)$
and $c_6=\d^{-1} \nu''(c) \max(\rho,2c)$, it is not hard to check that
for $(R_i,C_i) \in \Aa(\d,r,\rho,c)$ and $v \geq t \geq 0$,
\begin{eqnarray*}
\D I_1(v,t) &\le& c_3 \D C(v) \int_0^v e^{-\d (v-u)} du
 + c_4 \D R(v) \int_0^v e^{-\d |u-t|} du
\le c_5 [\D R(v) + \D C(v)]  \\
\D I_2(v,t) &\le&
\rho \nu''(c) \D C(v) \int_0^t e^{-\d (t-u)} du
+ c \nu''(c) \D R(t) \int_0^t e^{-\d |v-u|} du
\le c_6[\D R(v) + \D C(v) ]
\end{eqnarray*}
Thus, with $c_7 = c \rho K/3$ and $\b^2 \nu''(c) r^{-1} = r \leq 1/6$,
we have from \eqref{DL} that for $s \geq t \geq 0$,
\begin{eqnarray*}
\D\tilde C(s,t)&\le& \D\L (s,t)+ c_7 \int_t^s \D\L(s,v) dv
+ \b^2 \int_t^s e^{-\frac{1}{2}(s-v)} [\D I_1(v,t) + \D I_2(v,t) ] dv \\
&\le& \b^2 M_C [\D C(s)+\D R(s)] \,,
\end{eqnarray*}
where
$$
M_C := c_2 \sup_{\theta \geq 0} \theta e^{-\theta/2} +
c_7 c_2 \int_0^\infty u e^{-u/2} du + 2 (c_5+c_6) \,,
$$
and consequently, \eqref{eq:lipC} holds.

\medskip
Suppose $\d,r,\rho,c$ are such that
$\Psi$ is a contraction on $\Aa(\d,r,\rho,c)$, 
hence also on its non-empty subset $\Sa(\d,r,\rho,c)$. 
Starting at some $S_0=(R_0,C_0) \in \Sa(\d,r,\rho,c)$ 
consider the sequence $S_k=\Psi(S_{k-1})$, $k=1,2,\ldots$,
in $\Sa(\d,r,\rho,c)$. Since $\Psi$ is a contraction,
clearly $\{S_k\}$ is a Cauchy sequence for 
the uniform norm $\|\cdot\|$ of \eqref{eq:unorm}. Hence, 
$S_k \to S_\infty$ in the Banach space
$(\CC(\bD) \ts \CC (\R_+\ts\R_+),\|\cdot\|)$.
Note that $\Aa(\d,r,\rho,c)$ is a closed 
subset of this Banach space, so 
$S_\infty \in \Aa(\d,r,\rho,c)$. Further, fixing $\tau \geq 0$,
with $|(x,y)|:=|x|+|y|$, since $S_k \in \Sa(\d,r,\rho,c)$ we have that 
$$
\lim_{T \to \infty} \sup_{t,t' \geq T} 
|S_\infty(t+\tau,t)-S_\infty(t'+\tau,t')|
\leq 2 \|S_\infty-S_k\| + 
\lim_{T \to \infty} \sup_{t,t' \geq T} 
|S_k(t+\tau,t)-S_k(t'+\tau,t')| = 2 \|S_\infty-S_k\| \,.
$$
Taking $k \to \infty$ we deduce
that $\{S_\infty(t+\tau,t)\}$ is a Cauchy function from 
$\R_+$ to $[0,\rho] \ts [0,c]$, hence          
$S_\infty(t+\tau,t)$ converges as $t \to \infty$. With this applying for
each $\tau \geq 0$, we see that $S_\infty \in \Sa(\d,r,\rho,c)$ and
further that $S_\infty$ is the unique fixed point of the contraction 
$\Psi$ on the metric space $(\Sa(\d,r,\rho,c),\|\cdot\|)$.
By our construction of $\Psi$, it follows that the
fixed point $S_\infty=(R,C)$ of $\Psi$
satisfies \eqref{eqRs}-\eqref{eqCs}, from which we conclude that the 
unique solution of the latter equations is in $\Sa(\d,r,\rho,c)$. As
noted before, this  
yields the existence of $R_{\fdt}(\tau)$ of \eqref{Rfdt-ex}
and $C_{\fdt}(\tau)$ of \eqref{Cfdt-ex}, such that 
$(R_{\fdt},C_{\fdt}) \in \Ba(\d,r,\rho,c)$.

In the course of proving Proposition \ref{Psi} we found that on 
$\Sa(\d,r,\rho,c)$ the mapping $\Psi$ induces a mapping
$\Psi_{\fdt}:(R_{\fdt},C_{\fdt}) \to (\tilde R_{\fdt},\tilde C_{\fdt})$
such that
\begin{eqnarray*}
\tilde R_{\fdt} (\tau) &=& \hh{\L}(\tau) 
\sum_{n\ge 0}\beta^{2n} \sum_{\s\in\NC_n}
\int_{0\le \theta_1\le \cdots\le \theta_{2n}\le \tau}
\prod_{i\in \cro(\s)} \nu''(C_{\fdt}(\theta_i-\theta_{\sigma(i)}))
\prod_{j=1}^{2n} d\theta_j \,,\\
\tilde C_{\fdt} (\tau) &=& \hh{\L}(\tau) +
\b^2 \int_0^{\tau} \hh{\L}(\tau-v) \hh{I}_1(v) dv
+ \b^2 \int_0^\tau \hh{\L}(\tau-v) \hh{I}_2(v) dv \,,
\end{eqnarray*}  
where $\hh{\L}(\tau) = e^{-\mu \tau}$ for $\mu=1/2+\b^2 \hh{I}_0$ 
of \eqref{eqKFDT}, while 
$\hh{I}_1(v)$
and $\hh{I}_2(v)$
of \eqref{eqI1l}-\eqref{eqI2l}
are the two integrals on the right-hand-side of \eqref{eqCFDT}. In particular,
$\tilde R_{\fdt}$ and $\tilde C_{\fdt}$ are differentiable on $\R_+$, and
by \cite{GM} we have that for $\tau \geq 0$, 
\begin{equation}\label{eqfPsiR}
\tilde R_{\fdt}' (\tau) = -\mu \tilde R_{\fdt} (\tau) 
+ \b^2 \int_0^\tau \tilde R_{\fdt} (\tau-\theta) \tilde R_{\fdt} (\theta) 
\nu''(C_{\fdt}(\theta)) d\theta \,,
\end{equation}
with $\tilde R_{\fdt}(0)=1$, while 
\begin{equation}\label{eqfPsiC}
\tilde C_{\fdt}' (\tau) =-\mu \tilde C_{\fdt} (\tau) + \b^2 \hh{I}_1(\tau) 
+ \b^2 \hh{I}_2(\tau) \,,
\end{equation}
with $\tilde C_{\fdt}(0)=1$. Since the solution $(R,C)$ of 
\eqref{eqRs}-\eqref{eqCs} is a fixed point of $\Psi$, 
the corresponding pair $(R_{\fdt},C_{\fdt})$ is a fixed point of
$\Psi_{\fdt}$, which by \eqref{eqfPsiR}--\eqref{eqfPsiC} 
satisfies the FDT equations \eqref{eqRFDT}-\eqref{eqKFDT} with $\IJ=0$.

Recall that $\Ba(\d,r,\rho,c)$ consists of pairs $(R,C)$ of 
functions that are uniformly bounded, exponentially decaying to zero at infinity, and 
of fixed values at zero. It is not hard to verify that when $\IJ=0$, any solution of
\eqref{eqRFDT}-\eqref{eqKFDT} in this space, is a fixed point of $\Psi_{\fdt}$.
Further, $\Psi_{\fdt}$ is a contraction on $\Ba(\d,r,\rho,c)$, 
equipped with the supremum norm, since following 
the very same arguments we used 
in proving the Lipschitz estimates \eqref{eq:lipR}-\eqref{eq:lipC}
of Proposition \ref{FDT1}, we find that
$\Delta \tilde R_{\fdt} (\tau) \leq M_R \b^2 [ \D R_{\fdt} (\tau) 
+ \D C_{\fdt} (\tau) ]$ and 
$\Delta \tilde C_{\fdt} (\tau) \leq M_C \b^2 [ \D R_{\fdt} (\tau) 
+ \D C_{\fdt} (\tau) ]$ for all $\tau<\infty$, 
now with $\D f(\tau)=\sup_{\theta \leq \tau} |f(\theta)|$.
This of course proves the uniqueness of the solution of
\eqref{eqRFDT}-\eqref{eqKFDT} in $\Ba(\d,r,\rho,c)$, in case $\IJ=0$, as claimed.
\hfill \qed

\section{Study of the FDT equations}\label{FDTeqsec}

In this section, we 
complete the proof of Theorem
\ref{FDT} by relating the solutions
of the equations \eqref{eqRFDT}-\eqref{eqKFDT}
with the solution of \eqref{FDTDb}.

Specifically, we first prove Proposition \ref{studyFDTeq}
about existence, uniqueness, limiting value and exponential 
convergence of the solution of \eqref{FDTDb}.
Fixing $\b<\b_c$, we then know that  
the unique solution of \eqref{FDTDb} 
for $b=1/2$ and $\phi(x)=b+2 \b^2 \nu'(x)$, exists,
is twice continuously differentiable, 
positive, and of negative derivative, such that 
both converge exponentially fast to zero when $s \to \infty$,
with positive exponent $\epsilon_\b = \frac{1}{2}-2\b^2\nu''(0)$. 
Further, for $\b < \b_1 \leq \b_c$ sufficiently small, 
by Proposition \ref{FDT1}, the pair of functions 
$(R(\tau+t,t),C(\tau+t,t))$ converges for $t \to \infty$ to the
unique solution $(R_{\fdt}(\tau),C_{\fdt}(\tau))$ of 
\eqref{eqRFDT}-\eqref{eqKFDT}, with $\IJ=0$. Hence, our next proposition 
(whose proof is provided at the end of the section),
completes the proof of the theorem by showing that
for any $\b$, the solution
$(C_{\fdt}, R_{\fdt})$ of \eqref{eqRFDT}-\eqref{eqKFDT}
can be expressed in terms of the solution of \eqref{FDTDb} for
$b=1/2$ and $\phi(x)=\gamma+2\b^2 \nu'(x)$, provided the constant
$\gamma$ is chosen accordingly. That is, with $\IJ_\gamma$ matching 
the constant
$\IJ$ of \eqref{eqKFDT} while $\phi(\cdot)$ satisfies 
the condition \eqref{eq:Dcond} of Proposition \ref{studyFDTeq}.
Indeed, if $\b<\b_c$, then for $\gamma=b$ we have 
that $D_\infty=0$, $\IJ_\gamma=0$ and further, the resulting solution $(R,C)$ 
has the exponential decay property of Proposition 
\ref{FDT1} provided $\b$ is small enough (for the condition 
$b+2 \b^2 \nu'(1) > 2 \sqrt{b 2 \b^2 \nu''(1)}$ of Proposition \ref{studyFDTeq}
to hold).
\begin{prop}\label{relationeq} 
Suppose $D(s)$ is 
a positive, twice continuously differentiable and decreasing 
solution of the equation \eqref{FDTDb} for $b>0$, 
$D_\infty \geq 0$ and $\phi(\cdot)$ which satisfy the conditions of
Proposition \ref{studyFDTeq}. Then, $R(s)=-b^{-1} D'(s) \geq 0$,
$C(s)=C(-s)=D(s) \geq 0$ and $\mu=\phi(1)>0$ are such that
$R(0)=C(0)=1$ and for any $\gamma \in \R$,
\begin{eqnarray}
R'(\tau)&\!\!\!\! =\!\!\!\!
& - \mu R(\tau) + b \int_0^\tau R(\tau-\theta)
R(\theta) \phi'(C(\theta))d\theta ,\label{eqRFDTb}\\
C'(\tau)&\!\!\!\! = \!\!\!\!
& - \mu C(\tau)
+ b \int_0^\infty C (\tau-\theta) R (\theta) \phi'(C (\theta)) d\theta
+ b \int_\tau^\infty [\phi(C(\theta))-\gamma] R(\theta-\tau)d\theta + 
\IJ_\gamma
\label{eqCFDTb}\\
\mu&= & b + b 
\int_0^\infty \npsi_\gamma (C(\theta)) R(\theta) d\theta + \IJ_\gamma
 \,, \label{eqKFDTb}
\end{eqnarray}
for $\npsi_\gamma (x)=x\phi'(x)+\phi(x)-\gamma$ and 
$\IJ_\gamma =\gamma-b+D_\infty(\phi(D_\infty)-\gamma)$
(recall also that $b=\phi(D_\infty)(1-D_\infty)$).
\end{prop}

\nn
{\bf Proof of Proposition \ref{studyFDTeq}.}
Let $I=[0,1-D_\infty]$ and
$M_\infty(\R_+,I)$ denote the
space of continuous $I$-valued functions on $\R_+$ which take the
value $1-D_\infty$ at time zero.
Set $\phsi(x)=\phi(x+D_\infty)-\phi(D_\infty)$ and
for any $E \in M_\infty(\R_+,I)$, let $k(s)=b \phsi'(E(s))$ and
\begin{equation}\label{eq:Hdef}
H_s(E)=1+\sum_{n\ge 1} \sum_{\s\in NC_n}
\int_{0\le t_1\cdots \le t_{2n}\le s}
\prod_{i \in \cro(\s)}  k(t_{\s(i)}-t_i) \prod_{i=1}^{2n} dt_i \,.
\end{equation}
Note that the continuously differentiable
$s \mapsto H_s(E)$ is the solution of
\begin{equation}\label{eq:Hform}
\frac{d H_s(E)}{ds} =b \int_0^s \phsi'(E(s-v)) H_{s-v} (E) H_{v} (E) dv,
\qquad H_0(E)=1 \,.
\end{equation}
For $\mu =\phi(1)$ let
$$
\Phi(E)(s) =\big( 1-D_\infty- b\int_0^s e^{-\mu v} H_v(E) dv
\big) \vee 0 \,.
$$
Since $b\phsi'$ is non-negative on the interval $I$,
it follows that $H_s(E)$ is non-negative and consequently,
$s \mapsto \Phi(E)(s)$ is non-increasing and
also belongs to $M_\infty(\R_+,I)$. Further,
since $\phsi'$ is by assumption Lipschitz continuous on the compact set $I$
it follows from \eqref{eq:Hdef} that for each $T<\infty$ there exists
$\kappa_1(T)<\infty$ such that for all 
$E, \bar E \in M_\infty(\R_+,I)$ and $u \leq T$,
$$|H_u(E)-H_u(\bar E)|\le \kappa_1(T) \int_0^u |E(v)-\bar E (v)| dv$$
With $y\mapsto y \vee 0$ Lipschitz continuous, this implies that for
some finite, non-decreasing 
$\kappa_2(T) \geq 1$, all $E, \bar E \in M_\infty(\R_+,I)$ and 
$s \leq T$,
$$
|\Phi(E)(s)-\Phi(\bar E)(s)|\le \kappa_2 (T)\int_0^s |E(v)-\bar E (v)| 
dv\,.
$$
Thus, $\Phi(\cdot)$ is a contraction on
$M_\infty(\R_+,I)$ equipped with the weighted $L_1$-norm  
$\| E \|_* := \int_0^\infty |E(s)| w(s) ds$
for weight function $w(s)=\exp(-2 \int_0^s \kappa_2(u) du)>0$.
Consequently, we deduce that $\Phi(E)=E$ has a unique solution 
in $M_\infty(\R_+,I)$, denoted hereafter by $E_*(s)$. 
Let $\sigma_*>0$ be the first $s>0$ where
$\Phi(E_*)(s)=0$.
Note that the function $s \mapsto \Phi(E_*)(s)$ is continuously differentiable
on the interval $[0,\sigma_*)$, in which case also
$$
E_*(s)=1-D_\infty- b \int_0^s e^{-\mu v} H_v(E_*) dv \,,
$$
is twice differentiable, with $E_*'(s)=-b e^{-\mu s} H_s(E_*)$.
Hence, applying \eqref{eq:Hform} for $H_s(E_*)$,
we have that for $s \in [0,\sigma_*)$,
$$
E_* ''(s)= -\mu E_* '(s) - \int_0^s \phsi'(E_*(s-v)) E_*'(s-v) E_*'(v) dv.
$$
Integrating this equation we find that 
$$
E_*'(s)=E_*'(0) -\int_0^s [\mu +\phsi(E_*(s-v))-\phsi(E_*(0))] E_*'(v) dv \,.
$$
Since $E_*'(0)=-b$, $E_*(0)=1-D_\infty$ and 
$\mu-\phsi(E_*(0))=\phi(D_\infty)$, it follows that on $[0,\sigma_*)$
\begin{equation}\label{eqE}
E_*'(s)=-b -\int_0^s \phi(E_*(s-v)+D_\infty) E_*'(v) dv \,, \qquad
\E_*(0) = 1 - D_\infty \,.
\end{equation}

We next show that $\sigma_*=\infty$, which in view of \eqref{eqE}
yields that 
$D(s)=E_*(s)+D_\infty$ is a twice differentiable, strictly 
decreasing $[D_\infty,1]$-valued solution of \eqref{FDTDb}.
To this end, note that $H_s(E_*) \geq 1$ since $E_* \in M_\infty(\R_+,I)$ 
(see \eqref{eq:Hdef}).
Consequently, for all $v < \sigma_*$ both
$E_*'(v) \leq -b e^{-\mu v} < 0$ and 
$$
- \phi(E_*(s-v)+D_\infty) E_*'(v) \geq -\phi(D_\infty) E_*'(v) \,.
$$
Thus, from \eqref{eqE} we have that for $s<\sigma_*$,
\begin{equation}\label{amir90}
-b e^{-\mu s} \geq E'_*(s) \geq -b -\phi(D_\infty) (E_*(s) - E_*(0))
= - \phi(D _\infty) E_*(s) \,,
\end{equation}
since $b=\phi(D_\infty) (1-D_\infty) =\phi(D_\infty) E_*(0)$ by 
the definition of $D_\infty$. If 
$\sigma_*<\infty$, then as $s \uparrow \sigma_*$ we have by the
continuity of $E_*(\cdot)$ that 
$E_*(s) \to E_*(\sigma_*) = \Phi(E_*)(\sigma_*)=0$, in contradiction
with \eqref{amir90}.

The uniqueness of the solution to \eqref{FDTDb} follows 
from the uniqueness of the preceding $E_*(\cdot)$. Indeed, 
if a continuously differentiable $[D_\infty,1]$-valued $D(s)$ solves
\eqref{FDTDb}, then $E(s)=D(s)-D_\infty \in M_\infty(\R_+,I)$ and 
since $b=\phi(D_\infty) E(0)$ we have that 
\begin{equation}\label{eqE2}
E'(s)=-\phi(D_\infty) E(s) -\int_0^s \phsi(E(s-v)) E'(v) dv \,, \qquad
E(0) = 1 - D_\infty \,.
\end{equation}
Hence, $E(s)$ is twice differentiable, with $E'(0)=-b$, 
and differentiating \eqref{eqE2} we get that
\begin{equation}\label{eqRnew}
E''(s) = - \mu E'(s) - \int_0^s \phsi'(E(s-v)) E'(s-v) E'(v) dv
\end{equation}
(using the fact that $\phsi(E(0))+\phi(D_\infty)=\mu$). 
Thus, $H_s(E)=-b^{-1} e^{-\mu s} E'(s)$ solves \eqref{eq:Hform}.
This means that for all $s \geq 0$,
\begin{equation}\label{amir85}
E(s) =1-D_\infty- b\int_0^s e^{-\mu v} H_v (E) dv \,.
\end{equation}
Since $E(s) \geq 0$, so is the right-hand-side of \eqref{amir85},
that is $E=\Phi(E)$. Since $\Phi(\cdot)$ has a unique fixed point
in $M_\infty(\R_+,I)$, the solution of \eqref{FDTDb} must also be unique.

Next, the monotone and bounded function $D$ 
converges to some $x\in [D_\infty,1]$. Fix $M=M_\epsilon$ such that
$D(v) \leq x + \epsilon$ for all $v \geq M$. Then, for $s \geq 2M$
we have that
\begin{eqnarray*} 
\int_0^s \phi(D(v)) D'(s-v) dv
&=&\int_0^M \phi(D(v)) D'(s-v) dv +\int_M^s \phi(D(v)) D'(s-v) dv \\
&\ge& - \phi(1) (D(s-M)- D(s)) - \phi(x+\epsilon) (1-D(s-M))\\
&\ge& - \phi(1) (x+\epsilon-x) -\phi(x+\epsilon) (1-x) \,.
\end{eqnarray*}
Hence, for all $s \geq 2M$, 
$$
D'(s) \leq -b + \phi(1) \epsilon +\phi(x+\epsilon) (1-x) \,.
$$
Since $D(s)$ is bounded below, it follows that  
$\phi(1) \epsilon +\phi(x+\epsilon) (1-x) \geq b$.
Taking $\epsilon \downarrow 0$ we see that 
$x=D_\infty$ (since $x\ge D_\infty$ and
$D_\infty$ is the largest $y \in [0,1]$ for 
which $\phi(y) (1-y) \geq b$).

Finally, recall 
that the function $E(s)=D(s)-D_\infty$ is strictly positive, 
monotone decreasing with $E'(s)=-b e^{-\mu s} H_s(E)$ for
$H_s(E)$ of \eqref{eq:Hdef}. Since there 
$k(s)=b \phsi'(E(s)) \leq b \phi'(1)$ (by the 
assumed convexity of $\phi(\cdot)$),
the same argument as in \eqref{eqtHbd1} yields that for all $s \geq 0$,
$$
H_s(E) \leq (2\pi)^{-1} 
\int_{-2}^2 e^{\sqrt{b \phi'(1)} s x} \sqrt{4-x^2} dx
\leq c_1 e^{2 \sqrt{b \phi'(1)} s} \,.
$$ 
Consequently, if $\phi(1) > 2 \sqrt{b \phi'(1)}$, then 
both $E'(s)$ and $E(s)$ converge exponentially to zero for $s \to \infty$.
\hfill\qed

\nn{\bf Remark.} We note in passing that 
$\phsi(x) \geq \phsi'(0) x$ for all $x \in [0,1]$ when 
$\phi$ is convex, in which case by the monotonicity of $E(s)$,
$$
J(s) := E(s)^{-1} \int_0^s [-E'(v)] \phsi(E(s-v)) dv 
\ge (x^{-1}\phsi)(E(s))[E(0)-E(s)]\ge \phsi'(0) [1-D_\infty - E(s)] \,.
$$
Further, recall that from \eqref{eqE2},
$$
\frac{d \log E}{ds} (s) = \frac{E'}{E} (s) = J(s)-\phi(D_\infty) \,,
$$ 
hence in case $\phi(D_\infty) = \phi'(D_\infty) (1-D_\infty)$, we find that 
$E(s) \geq E(0) \exp(-\phi'(D_\infty) \int_0^s E(u) du)$, which thus 
does not converge to zero exponentially fast as $s \to \infty$.

\medskip
\nn {\bf Proof of Proposition \ref{relationeq}.} 
Since $D(0)=1$ and $D(\theta) \downarrow D_\infty$, 
for $R=-b^{-1} D'$ and any $\gamma \in \R$ we have that 
$$
\gamma b \int_0^\infty R(\theta) d\theta = \gamma (1-D_\infty) = \IJ_\gamma
-\IJ_0 \,,
$$
hence we may and shall assume hereafter that $\gamma=0$. Further, 
it is easy to 
check that for $\npsi_0(x)=[x\phi(x)]'$ 
$$
\int_0^\infty \npsi_0 (D(\theta)) D'(\theta) d\theta = 
D_\infty \phi(D_\infty) - D(0) \phi(D(0)) \,, 
$$
so the choices of $\mu=\phi(1)=\phi(D(0))$ and 
$\IJ_0=-b+D_\infty \phi(D_\infty)$ guarantee
that $C=D$ and $R=-b^{-1} D'$ satisfy \eqref{eqKFDTb}. 
Recall that while proving Proposition \ref{studyFDTeq}
we have seen that $D'(0)=-b$ and for all $\tau \geq 0$
$$
D''(\tau)=-\mu D'(\tau) -\int_0^\tau D'(\tau-\theta) D'(\theta) 
\phi'(D(\theta)) d\theta 
$$
(c.f. \eqref{eqRnew}), from which it immediately follows that
$C=D$ and $R=-b^{-1} D'$ satisfy \eqref{eqRFDTb}.

Turning to the remaining task of verifying that \eqref{eqCFDTb}
holds, note that since $D(\cdot)$ is bounded and converges to 
$D_\infty$, 
$$
\int_\tau^\infty D(\theta-\tau) D'(\theta) \phi'(D(\theta)) d\theta
+ \int_\tau^\infty \phi(D(\theta)) D'(\theta-\tau)d\theta 
= \phi(D_\infty) D_\infty - \phi(D(\tau)) D(0) \,.
$$
Similarly, 
$$
\int_0^\tau D(\tau-\theta) D'(\theta) \phi'(D(\theta)) d\theta
- \int_0^\tau \phi(D(\theta)) D'(\tau-\theta)d\theta 
= \phi(D(\tau)) D(0) - \phi(D(0)) D(\tau) \,.
$$
Hence, using the symmetry $C(s)=C(-s)$, upon substituting 
in \eqref{eqCFDTb} $C=D$, $R=-b^{-1} D'$ and  
the preceding two identities, it is not hard to  
verify that \eqref{eqCFDTb} holds for $\gamma=0$ if 
$$
D'(\tau)=-\phi(D_\infty) D_\infty + 
\IJ_0  
- \int_0^\tau \phi(D(\theta)) D'(\tau-\theta)d\theta \,,
$$
which in view of \eqref{FDTDb} is merely the statement that
$b = \phi(D_\infty) D_\infty - \IJ_0$. 
Our choice of $\IJ_0$ guarantees that the latter identity applies,
thus completing the proof of the proposition.
\hfill\qed


\section{About the FDT regime at all temperatures}\label{sec-heur}
We return to the equations \eqref{eqCs}--\eqref{eqZs}
and explain why \eqref{FDTDb} is the natural candidate
for describing the function $C_{\fdt}(\cdot)$ for all $\b$.
To this end, set for $s \geq t$,
$$G(s,t)=R(s,t)-2\partial_t C(s,t)$$
and
$$I(s,t)=\beta^2\int_0^t [
C(t,u) G(s,u) \nu''(C(s,u)) + \nu'(C(s,u)) G(t,u)
] du
-2\beta^2\nu'(C(s,0))C(t,0).$$
Since $\psi(x)=[x\nu'(x)]'$ and $C(s,s)=1$, it follows that
\begin{eqnarray}\label{eq:Iss}
I(s,s)&=&\b^2 \int_0^s \psi(C(s,u)) G(s,u) du -
2\beta^2\nu'(C(s,0))C(s,0) \nonumber \\
&=& \b^2 \int_0^s \psi(C(s,u)) R(s,u) du - 2\b^2 \nu'(1) = \mu - \frac{1}{2}
- 2\b^2 \nu'(1) \,,
\end{eqnarray}
for $\mu$ of \eqref{eqZs}. Further, by similar reasoning, 
\begin{eqnarray*}
2\b^2 \nu'(C(s,t)) + I(s,t) &=&
2 \b^2 \int_0^t \partial_u [C(t,u) \nu'(C(s,u))] du 
+ 2\beta^2\nu'(C(s,0))C(t,0)  + I(s,t) \\
&=& \b^2 \int_0^t [C(t,u) R(s,u) \nu''(C(s,u)) + \nu'(C(s,u)) R(t,u)] du \,.
\end{eqnarray*}
Thus, in these notations \eqref{eqRs}-\eqref{eqZs} are equivalent to 
having for $s \geq t$, 
\begin{eqnarray}
\partial_s R(s,t)&=&-[\rho+I(s,s)] R(s,t) 
+\beta^2\int_t^s R(s,u)R(u,t) \nu''(C(s,u)) du \,,\label{eqRs2}\\
\partial_s C(s,t)&=&-[\rho+I(s,s)] C(s,t) 
+\beta^2\int_t^s C(u,t) R(s,u) \nu''(C(s,u)) du 
+2\beta^2 \nu'(C(s,t)) +I(s,t)\label{eqCs2}
\,,
\end{eqnarray}
with $R(t,t)=C(t,t)=1$ and $\rho=2^{-1} +2\beta^2 \nu'(1)$.

We note in passing that $G(s,t)$ is such that for $s \geq t$,
\begin{eqnarray}\label{eqGn}
\partial_s G(s,t) &=& - 2 \partial_t I(s,t) - 
[\rho+I(s,s)] G(s,t) + 2 \beta^2 \nu''(C(s,t)) G(s,t) \\
&+& \beta^2 \int_t^s G(u,t) [ G(s,u)\nu''(C(s,u) 
+ 2 \partial_u C(s,u) \nu''(C(s,u)) ] du \,. \nonumber
\end{eqnarray}
The physics prediction is that 
$G(t+\tau,t) \to 0$ as $t \to \infty$ while $\tau$ 
is fixed, for any finite $\beta$
(this is the famous FDT relation). As a result, comparing 
$I(s,t)$ and $I(t,t)$, we 
further expect that $I(t+\tau,t) \to \wIJ$ when $t \to \infty$ while 
$\tau$ is fixed. We next show that the latter ansatz 
results with the existence of an FDT solution 
$(R_{\fdt},C_{\fdt})$ such that $R_{\fdt}=-2 C_{\fdt}'$ and
$C_{\fdt}$ solves \eqref{FDTDb} for 
$b=1/2$ and $\phi(x)=\gamma+2\b^2 \nu'(x)$, where 
$\gamma=\wIJ+1/2$.
\begin{prop}\label{ansatz} 
Assume that given the continuous function $I(s,t)$ 
there exists a continuously differentiable  
solution $(R_I(s,t),C_I(s,t))$, $s \geq t \geq 0$,
of \eqref{eqRs2}--\eqref{eqCs2} 
with $\rho=2^{-1} +2\beta^2 \nu'(1)$, the 
initial conditions $R(t,t)=C(t,t)=1$ and uniformly bounded
$C(s,t)$. Further, suppose that  
for any $T<\infty$ 
\begin{equation}\label{Iass}
\lim_{t \to \infty} \sup_{\tau \in [0,T]} \, |I(t+\tau,t) - \wIJ| = 0 \,,
\end{equation} 
where the constant $\wIJ$ is such that 
\begin{equation}\label{eq:DcondI}
\sup_{0 \leq  x \leq 1} \, 
\{ (\wIJ+\frac{1}{2}+2\beta^2 \nu'(x))(1-x) \} \ge \frac{1}{2} \,.
\end{equation} 
Then, 
$C_I(t+\tau,t) \to C_{\fdt}(\tau)$ as $t \to \infty$, 
uniformly in $\tau \in [0,T]$, where $C_{\fdt}(\tau)$ is the unique solution 
of \eqref{FDTDb} for $b=1/2$ and $\phi(x)=\gamma+2\b^2 \nu'(x)$,
with $\gamma=\wIJ+1/2$. Further, 
$R_I(t+\tau,t) \to R_{\fdt}(\tau) = -2 C_{\fdt}'(\tau)$ as 
$t \to \infty$, uniformly in $\tau\in [0,T]$.
\end{prop}
\nn{\bf Remark.} In Theorem \ref{FDT} we
circumvent the difficulty of showing the ansatz that 
$I(t+\tau,t)$ converges to a ($\b$-dependent) 
constant at any value of $\b$, by verifying 
that at very high temperature, i.e. sufficiently small $\b$,
the exponential decay to zero in $s-t$ of $(R,C)$ 
results with the convergence to zero of $I$. Proposition 
\ref{ansatz} shows that this is the only obstacle to 
extending our results about the FDT regime to all temperatures.

Before proving the proposition, consider its consistency with 
our choice of 
$\gamma(\b)$ for $\b \geq \b_c$ based on $q(\b)$ of \eqref{qEA}.
Specifically, considering in \eqref{eq:Iss} 
the contribution to the integral from $u \in [s-M,s]$ with $M<\infty$ 
arbitrarily large, we expect the FDT solution $(R_{\fdt},C_{\fdt})$ to
contribute $-2\b^2 q(\b) \nu'(q(\b))$ to the limiting constant 
$\wIJ=\gamma(\b)-1/2$ (for which \eqref{eq:DcondI} holds by our
choice of $\gamma(\b)$).  
We further expect $\wIJ$ to be the sum of this FDT contribution 
and a non-negative contribution from the aging regime (i.e.
the integral over $u \in [0,s-M]$). Given the 
relation between $\mu$ and $\IJ_\gamma$ in 
Proposition \ref{relationeq}, we deduce that 
$\IJ_\gamma$ is exactly the contribution of the aging regime 
to $\wIJ$. In conclusion, we should 
have $\IJ_\gamma=0$ for $\b<\b_c$, $\IJ_\gamma \geq 0$
at $\b=\b_c$ and $\IJ_\gamma>0$ for $\b>\b_c$.
This is indeed the case, for when $\b<\b_c$ we have
$\gamma=b=1/2$ and $D_\infty=0$ leading to $\IJ_\gamma=0$ while
for $\b \geq \b_c$, since 
$\nu'(x) \leq x \nu''(x)$ for $x \geq 0$ with strict inequality when 
$x>0$, our choice of $\gamma(\b)$ leads to 
$$
2 I_\gamma=4\b^2 (1-q(\b)) [\nu''(q(\b))-\nu'(q(\b))] - 1 
\geq 4 \b^2 (1-q(\b))^2 \nu''(q(\b)) -1 \geq 0 \,,
$$
with a strict inequality whenever $q(\b)>0$ (in particular, for all
$\b>\b_c$).

\nn {\bf Proof of Proposition \ref{ansatz}.} 
Fixing hereafter the values of $\b$ and $\rho=1/2+2\b^2\nu'(1)$ and 
the initial conditions $R(t,t)=C(t,t)=1$, 
a re-run of the argument at the end of the proof of 
Proposition \ref{prop-sphere} shows that per given 
continuous function $I(s,t)$, the 
system of equations \eqref{eqRs2}--\eqref{eqCs2} admits 
at most one bounded solution, denoted $(R_I(s,t),C_I(s,t))$ on any
compact interval $0 \leq t \leq s \leq T$. Consequently, 
it has at most one continuous solution $(R_I,C_I)$ for all
$s \geq t \geq 0$. 

In particular, in case $I(s,t)=\wIJ$ is 
a constant that satisfies \eqref{eq:DcondI}, we know 
that \eqref{eq:Dcond} holds for
$b=1/2$, $\gamma=\wIJ+1/2$ and $\phi(x)=\gamma+2\b^2 \nu'(x)$. 
With $D(\tau)$ denoting the unique solution of \eqref{FDTDb}
for these parameters, we claim that 
$R_{\wIJ}(s,t)=-b^{-1} D'(s-t)$ and $C_{\wIJ}(s,t)=D(s-t)$ 
is then a solution of \eqref{eqRs2}--\eqref{eqCs2} (hence its unique
solution). Indeed, recall that $D(0)=1$, $D'(0)=-b$ so the given initial 
conditions for $(R_{\wIJ},C_{\wIJ})$ hold. Further,
these settings result with $\mu:=\rho+\wIJ=\phi(1)$, so 
upon taking $u=t+\theta$ and $s=t+\tau$ we see that 
the equation \eqref{eqRs2}, which is a direct 
consequence of \eqref{eqRFDTb} of Proposition \ref{relationeq},
clearly holds. By the same transformations
we see that our proposed solution satisfies \eqref{eqCs2} provided that
\begin{equation}\label{FDTDc}
D'(\tau)=-\phi(1) D(\tau) - 
\int_0^\tau D(\tau-\theta) D'(\theta) \phi'(D(\theta)) d\theta
+ \phi(D(\tau)) - b \,.
\end{equation}
Noting that by integration by parts, 
$$
\int_0^\tau D(\tau-\theta) D'(\theta) \phi'(D(\theta)) d\theta
= - D(\tau-\theta) \phi(D(\theta)) \vert_0^\tau 
+\int_0^\tau D'(\tau-\theta) \phi (D(\theta)) d\theta \,,
$$
and having $D(0)=1$, we see that \eqref{FDTDc} is equivalent to 
\eqref{FDTDb} hence holds as well.

Fixing $T<\infty$, for each $\xi$ positive let
$\bD_\xi:=\{(s,t): \xi \leq t \leq s \leq t + T \} \subset \R_+ \ts \R_+$.
We proceed by showing that per fixed $\xi \geq 0$ 
and $T<\infty$ the mapping $I \mapsto (R_I,C_I)$ is
Lipschitz with respect to the supremum norm over $\bD_\xi$
with a Lipschitz constant that is independent of $\xi$.
More precisely, given 
two bounded functions $I(s,t)$ and $\bar I(s,t)$ 
to which correspond continuously differentiable solutions
$(R_I,C_I)$ and $(R_{\bar I},C_{\bar I})$ of 
\eqref{eqRs2}--\eqref{eqCs2}, such that 
$C_I$ and $C_{\bar I}$ are uniformly bounded,
we let
$\Delta_\xi I := \sup \{|I(s,t)- \bar I(s,t)| : (s,t) \in \bD_\xi\}$,
and show that 
\begin{eqnarray}\label{eq:RIbd}
\Delta_\xi R &:=& \sup \{
|R_I(s,t)- R_{\bar I}(s,t)| : (s,t) \in \bD_\xi \}
 \leq \kappa_0 \Delta_\xi I ,\\
\Delta_\xi C &:=& \sup \{ |C_I(s,t)- C_{\bar I}(s,t)| 
: (s,t) \in \bD_\xi \} \leq \kappa_0 \Delta_\xi I,
\label{eq:CIbd}
\end{eqnarray}
where $\kappa_0$ is a finite constant 
that depends only on $T$, $\b$, $\nu'$, $\rho$ and
the uniform bound on $I$, $\bar I$, $C_I$ and $C_{\bar I}$
for $(s,t) \in \bD_\xi$. Considering
the bounded ${\bar I}=\wIJ$, we have seen already that 
$R_{\bar I}(t+\tau,t)=R_{\fdt}(\tau)$ and
$C_{\bar I}(t+\tau,t)=C_{\fdt}(\tau)$ 
with $C_{\bar I}$ bounded.
Our assumption \eqref{Iass} then amounts to 
$\Delta_\xi I \to 0$ as $\xi \to \infty$. 
In particular, this implies that
$I$ is also uniformly bounded on $\bD_\xi$, so 
from \eqref{eq:RIbd} and \eqref{eq:CIbd} we get
that $\Delta_\xi R \to 0$ and $\Delta_\xi C \to 0$, 
which are easily seen to match our desired conclusion.

We thus complete the proof by verifying the bounds of
\eqref{eq:RIbd} and \eqref{eq:CIbd}, by an argument 
similar to the one we used for proving uniqueness 
of the system of equations \eqref{eqRs}-\eqref{eqZs}.
To this end, note first that $R_I=\Lambda_I H_{C_I}$ for  
$\Lambda_I(s,t)=\exp(-\int_t^s(\rho+I(u,u)) du)$
and $H_C(s,t)$ of \eqref{gmfor}. Thus the assumed 
uniform bounds for $I$ and $C_I$ on $\bD_\xi$
imply that $R_I$ is also uniformly bounded on 
$\bD_\xi$ by a constant that depends only on 
$T$, $\b$, $\rho$, $\nu''$ and the corresponding 
uniform bounds for $I$, $C_I$.
Of course, the same applies for $R_{\bar I}$.
Next let 
$\Delta I(v,u)=|I(v,u)-{\bar I}(v,u)|$,
$\Delta R(v,u)=|R_I(v,u)-R_{\bar I}(v,u)|$,
$\Delta C(v,u)=|C_I(v,u)-C_{\bar I}(v,u)|$ and for
$(s,t) \in \bD_\xi$ set
$$
h(s,t)=\int_t^s [ \Delta R(s,u) + \Delta C(s,u) ] du \,.
$$
Then, similarly to the derivation of \eqref{tutuR} and \eqref{tutuC},
upon considering the difference between the integrated form of 
\eqref{eqRs2}--\eqref{eqCs2} for our solutions 
$(C_I,R_I)$ and $(C_{\bar I},R_{\bar I})$, we find that
for any $(s,t) \in \bD_\xi$
\begin{eqnarray}
\Delta R(s,t)&\le& \kappa_1 
[\int_t^s \Delta R(v,t) dv +\int_t^s h(v,t) dv + \int_t^s \Delta I(v,v) dv]
\,,
\label{tutuR2}\\
\Delta C(s,t)&\le& \kappa_1 [
\int_t^s \Delta C(v,t) dv +\int_t^s h(v,t) dv  + \int_t^s \Delta I(v,t) dv
+\int_t^s \Delta I(v,v) dv]
\,,
\label{tutuC2}
\end{eqnarray}
where the positive 
$\kappa_1 < \infty$ depends only on $T$, $\b$, $\rho$, $\nu'(\cdot)$ and
the maximum of $I$, 
$|R_I|$, $|C_I|$, $|\bar I|$, $|R_{\bar I}|$ and $|C_{\bar I}|$ on $\bD_\xi$.
Replacing $t$ by $u$, summing these two inequalities and
then integrating the result over $u \in [t,s]$ yields that 
\begin{equation}\label{hIbd}
h(s,t) \leq  \kappa_2 [
\int_t^s h(v,t) dt + \int_t^s \Delta I(v,v) dv + \int_t^s \int_t^v \Delta I(v,u)du dv ] \,,
\end{equation}
for any $(s,t) \in \bD_\xi$,
with a finite, positive constant $\kappa_2$ (of the same type of 
dependence as $\kappa_1$).
Since $h(t,t)=0$ for all $t$, 
we get from \eqref{hIbd} by Gronwall's lemma that 
$h(s,t) \leq \kappa_3 \Delta_\xi I$ for 
some finite, positive 
$\kappa_3$ (of same dependence type as $\kappa_1$) and
all $(s,t) \in \bD_\xi$. Recall that 
$\Delta R(t,t)=\Delta C(t,t)=0$ due to the given initial 
conditions, so upon plugging into \eqref{tutuR2} and \eqref{tutuC2} 
our uniform bound on $h(s,t)$, 
we complete the proof of \eqref{eq:RIbd}--\eqref{eq:CIbd} by
yet another application of Gronwall's lemma.
\hfill \qed

\end{document}